\documentclass{article}
\usepackage[utf8]{inputenc}
\usepackage{amsmath,amsthm,titling,enumerate,amssymb,etoolbox,mathrsfs,graphicx,xcolor,cancel,harpoon,tikz-cd,hyperref,mathtools}

\hypersetup{
    colorlinks=true,
    linkcolor=blue}

\title{$L^2$-estimates for the Dirac-Dolbeault operator and Bergman kernel asymptotics on some classes of non-compact complex manifolds}
\author{Ming-Yuan Chang}

\usepackage[margin=2.5cm]{geometry}

\def\N{\operatorname{\mathbb{N}}} 
\def\Z{\operatorname{\mathbb{Z}}}

\def\R{\operatorname{\mathbb{R}}}
\def\C{\operatorname{\mathbb{C}}}

\def\im{\operatorname{Im}}
\def\ker{\operatorname{Ker}}

\def\id{\operatorname{id}}

\def\tr{\operatorname{Tr}}

\def\Dom{\operatorname{Dom}}
\def\Spec{\operatorname{Spec}}

\newcommand{\qtext}[1]{\quad\text{#1}\quad}

\newcommand{\matrixtwo}[4]{ \left( \begin{array}{cc}
#1 & #2  \\
#3 & #4 
\end{array} \right)}

\newcommand{\vectortwo}[2]{ \left( \begin{array}{c}
#1   \\
#2 
\end{array} \right)}

\newtheoremstyle{mystyle}
  {6pt}{15pt}
  {}
  {}
  {\bf}
  {.}
  {1em}
  {}
  
\theoremstyle{mystyle}
\newtheorem{thm}{Theorem}[section]  
\newtheorem{definition}[thm]{Definition}
\newtheorem{lemma}[thm]{Lemma}
\newtheorem{proposition}[thm]{Proposition}
\newtheorem{corollary}[thm]{Corollary}
\newtheorem{remark}[thm]{Remark}

\newtheorem{example}[thm]{Example}

\newtheorem{assumption}[thm]{Assumption}

\date{}
\AtEndEnvironment{sol}{\null\hfill\qedsymbol}
\AtEndEnvironment{pf}{\null\hfill\qedsymbol}

\everymath{\displaystyle}

\begin{document}
\maketitle
\begin{abstract}
For high power $k$, the $L^2$-estimates for the Dirac-Dolbeault operator with coefficient $L^k\otimes E$ can be obtained from the Bochner-Kodaira-Nakano identity if $L$ has positive curvature. In this article, we generalize the classical method to obtain $L^2$-estimates for mixed curvature case, and give a bound to the extra error term. Modifying the $L^2$-estimates and existence theorem for $\bar{\partial}$-operator, we can get a local spectral gap of Kodaira Laplacian $\Box$ and thus a full asymptotic expansion for Bergman kernel.
\end{abstract}
\tableofcontents

\section{Introduction}

The asymptotic behavior of Bergman kernel with coefficient $L^k$ for a line bundle $L$ and $k$ large has important consequences in complex geometry and geometric quantization, such as the classical Kodaira embedding theorem and properties of Toeplitz operators. We refer to \hyperlink{ref:MZ}{[MZ]} \hyperlink{ref:HH}{[HH]} \hyperlink{ref:HMM}{[HMM]}.

In \hyperlink{ref:DLM}{[DLM]}, Dai-Liu-Ma established the full off-diagonal expansion for the Bergman kernel of the $\text{spin}^c$ Dirac operator associated to high powers of an ample line bundle over symplectic manifolds and orbifolds. In \hyperlink{ref:MM2}{[MM2]}, Ma-Marinescu continued to established the case that the line bundle curvature is nondegenerate with mixed curvature, that is, with negative and positive eigenvalues. Their method starts with deriving the spectral gap for the $\text{spin}^c$ Dirac operator on compact manifolds, then use a sequence of Sobolev estimate for the resolvent of the rescaled operator to analyze the residue formula. The method works generally for larger class of non-compact manifold if there are some uniform conditions on torsion and curvature to obtain global spectral gap \hyperlink{ref:MM1}{[MM1]}.

In \hyperlink{ref:HM}{[HM]}, Hsiao-Marinescu used PDE methods to obtain full diagonal expansion for general Hermitian manifolds under the assumption of small local spectral gap. For line bundles with semi-positive curvature, the classical Hörmander's $L^2$-estimates for $\bar{\partial}$-operator \hyperlink{ref:Hör}{[Hör]} gives such local spectral gap for complete Kähler manifolds and weakly pseuodoconvex domains.

A natural question is to find $L^2$-estimates and geometric conditions to obtain local spectral gap for mixed curvature case, or even allowing degeneracy. 

It is well-known that in holomorphic case, the estimate for semi-positive curvature can be derived using the Bochner-Kodaira-Nakano formula after twisting canonical bundle:
\[\frac{3}{2}(||\bar{\partial}s||^2+||\bar{\partial}^*s||^2)\ge\langle [\sqrt{-1}R^{E\otimes K_M^*},\Lambda]\tilde{s},\tilde{s}\rangle-\frac{1}{2}(||\mathcal{T}^*\tilde{s}||^2+||\bar{\mathcal{T}}\tilde{s}||^2+||\bar{\mathcal{T}}^*\tilde{s}||^2)\]
where $\tilde{\cdot}$ is the natural map from $\bigwedge\nolimits^{0,q}\otimes E$ to $\bigwedge\nolimits^{n,q}\otimes E\otimes K_M^*$.

As mentioned in \hyperlink{ref:MM2}{[MM2, Rmk.1.6]}, by using it the estimate for a fixed metric seems out of reach, and they used other formulas to deduce spectral gap.

However, for $R^L$ nondegenerate of signature $(d-q,q)$, pointwise if we unitarily diagonalize $R^L=\sqrt{-1}\sum_i\lambda_i w^i\wedge\bar{w}^i$ and $\lambda_1\ge\lambda_2\ge\cdots\ge\lambda_d$, then
\[[\sqrt{-1}R^L,\Lambda]w^I\wedge \bar{w}^J=(\sum_{i\in I}\lambda_i+\sum_{j\in J}\lambda_j-\sum_{l=1}^d\lambda_d)w^I\wedge\bar{w}^J\]
So if we pick $I$ to be the index set of all positive eigenvalue, the operator becomes positive definite on all $|J|\neq q$. The corresponding operation is twisting by the bundle $F_{p}$ generated by $w^I$ at every point. The problem is that unlike $K_M$, $F_{p}$ is not a holomorphic bundle, so we choose to use the bundle $\bigwedge\nolimits_{p,0}$ of $(p,0)$-vectors as coefficients, and there exists an error term $|\bar{\partial}\Psi^{F_p}|^2$ compared to the usual formula.

This article calculates the error term $|\bar{\partial}\Psi^{F_p}|^2$ and gives a bound to it. The following is the main result, which is a generalization of the usual estimate using Bochner-Kodaira-Nakano formula.
\begin{thm}
Let $(M,\Theta)$ be a Hermitian manifold of dimension $d$. $E$ and $L$ are holomorphic vector bundle and line bundle with Hermitian metrics. Assume $\lambda_p>\lambda_{p+1}$ and $[\sqrt{-1}R^{E\otimes\bigwedge\nolimits_{p,0}},\Lambda](z)\ge-\lambda^E_p(z)\id^{E\otimes \bigwedge\nolimits_{p,0}}(z)$ on $\bigwedge\nolimits^{p,q}$, then we have the following estimate for Dirac-Dobeault operator $D_k$:
\[\frac{1}{2}||D_ks||^2\ge \int_M\prescript{}{k}{f}_{p,q}|s|^2\]
for all $s\in\Omega^{0,q}_c(M,E\otimes L^k)$, where 
\[\prescript{}{k}{f}_{p,q}=
\begin{cases}
k\lambda_{p,q}-\lambda^E_p-|\bar{\partial}\Psi^{F_p}|^2 &\text{ if $M$ is Kähler}\\
\frac{2}{3}k\lambda_{p,q}-\frac{2}{3}\lambda^E_p-C_d|\partial\Theta|^2-|\bar{\partial}\Psi^{F_p}|^2 &\text{ if $M$ is non-Kähler}
\end{cases}\]
$C_d$ is a constant only depending on dimension $d$, and $\lambda_{p,q}=\sum_{i=1}^p\lambda_i-\sum_{j=1}^{d-q}\lambda_j$. 

Moreover, let $R:T^{1,0}\to T^{1,0}$ be the type change of $R^L$, we have the following bound for the error term $|\bar{\partial}\Psi^{F_p}|^2$:
\[|\bar{\partial}\Psi^{F_p}|^2\le(\lambda_{p+1}-\lambda_p)^{-2}|\bar{\partial}R|^2\]
\end{thm}
For $R^L$ nondegenerate with signature $(q,d-q)$, take $p=d-q$, and under appropriate geometric assumption we can get a local spectral gap of growth order $k$, and deduce the existence of full off-diagonal expansion.

In section \ref{section:preliminary}, we will derive basic equations in complex geometry.

Section \ref{section:linear algebra} consists of the construction of the twisting section, its property and the calculation for the corresponding error term.

Main result is proved in section \ref{section:L^2 estimate}. We will also show that the usual proof of $L^2$ estimate and existence theorem for $\bar{\partial}$-operator by Hörmander \hyperlink{ref:Hör}{[Hör]} can be modified to obtain that of the Dirac-Dolbeault operator $D$, under a general assumption (\ref{domain assumption}) for $\Dom\bar{\partial}$ and $\Dom\bar{\partial}^*$.

In section \ref{section:Bergman kernel}, we discuss the geometric conditions so that the main theorem applies and (\ref{domain assumption}) is satisfied. For manifolds with boundary we modify the calculation and obtain the corresponding boundary term with error terms. Finally, using the localization process in \hyperlink{ref:MM1}{[MM1]}, we show that asymptotics under local spectral gap condition is equivalent to the model case (\ref{model case}), which has a global spectral gap, and deduce a full off-diagonal expansion.


\section{Preliminaries: Bochner–Kodaira–Nakano identity}\label{section:preliminary}
The purpose of this section is to fix the notation and to derive some well-known geometric equations, especially the Bochner Kodaira Nakano inequality.

\[\mathscr{E}_{\eta}\tau=\eta\wedge\tau\]
\[\langle\mathscr{I}_{\eta}\tau_1,\tau_2\rangle=\langle\tau_1,\mathscr{E}_{\bar{\eta}}\tau_2\rangle\qquad\mathscr{I}_{\eta}=\mathscr{E}^*_{\bar{\eta}}\]
\[\mathscr{I}_{d\bar{z}^j}dz^I=g^{i\bar{j}}\delta_{i,K}^{I}\,dz^K\qquad \mathscr{I}_{dz^i}d\bar{z}^J=g^{i\bar{j}}\delta_{\bar{j},\bar{L}}^{\bar{J}}\,d\bar{z}^L\]
Define $\mathscr{I}_{X}$ to be the unique operator on forms satisfying graded Leibniz rule and for $1$-form $\eta$,
\[\mathscr{I}_X(\eta)=(\eta,X)\]
Then we have the relation
\[g^{p\bar{q}}\mathscr{I}_{\partial_p}=\mathscr{I}_{d\bar{z}^q}\qquad g^{p\bar{q}}\mathscr{I}_{\partial_{\bar{q}}}=\mathscr{I}_{dz^p}\]
Also define
\[\langle\mathscr{E}_{X}\tau_1,\tau_2\rangle=\langle\tau_1,\mathscr{I}_{\bar{X}}\tau_2\rangle\]
\[g^{p\bar{q}}\mathscr{E}_{\partial_p}=\mathscr{E}_{d\bar{z}^q}\qquad g^{p\bar{q}}\mathscr{E}_{\partial_{\bar{q}}}=\mathscr{E}_{dz^p}\]
\[d g_{I\bar{J}}=\sum_{i\in I,j\in J}dg_{i\bar{j}}(\operatorname{adj}g_{I\bar{J}})^{i\bar{j}}=\bar{\partial}g_{i\bar{j}}\delta^{iK}_{I}\delta^{\bar{j}\bar{L}}_{\bar{J}}g_{K\bar{L}}\]
If $|I|=|J|=p$, and at some point $g_{i\bar{j}}=\delta_{ij}$, then
\[d g_{I\bar{J}}=\begin{cases}
0&,|I\cap J|<p-1\\
\sum_{i\in I}dg_{i\bar{i}}&,I=J\\
\sum_{i\in I,j\not\in J}\delta^{i,I-i}_I\delta^{\bar{j},\overline{J-j}}_{\bar{J}}dg_{i\bar{j}}&,|I\cap J|=p-1
\end{cases}\]
We define the supercommutator of two $\Z_2$-graded operator:
\[[A,B]=AB-(-1)^{|A||B|}BA=(-1)^{|A||B|+1}[B,A]\]
\[[\mathscr{E}_{\eta_1},\mathscr{E}_{\eta_2}]=0=[\mathscr{I}_{\eta_1},\mathscr{I}_{\eta_2}]\]
For $1$-forms $\eta_1,\eta_2$,
\[[\mathscr{E}_{\eta_1},\mathscr{I}_{\bar{\eta}_2}]=\langle\eta_1,\eta_2\rangle\]
In particular
\[[\mathscr{E}_{dz^i},\mathscr{I}_{dz^k}]=0=[\mathscr{E}_{d\bar{z}^j},\mathscr{I}_{d\bar{z}^l}]\]
\[[\mathscr{E}_{dz^i},\mathscr{I}_{d\bar{z}^j}]=g^{i\bar{j}}=[\mathscr{E}_{d\bar{z}^j},\mathscr{I}_{dz^i}]\]
\[L=\mathscr{E}_{\sqrt{-1}g_{i\bar{j}}\,dz^i\wedge d\bar{z}^j}=\sqrt{-1}g_{i\bar{j}}\mathscr{E}_{dz^i}\mathscr{E}_{d\bar{z}^j}\qquad \Lambda=L^*=-\sqrt{-1}g_{i\bar{j}}\mathscr{I}_{dz^i}\mathscr{I}_{d\bar{z}^j}\]
\[[\mathscr{E}_{\eta},L]=0\]
\[[\mathscr{I}_{dz^i},\mathscr{E}_{dz^k}\mathscr{E}_{d\bar{z}^l}]=-g^{i\bar{l}}\mathscr{E}_{dz^k}\qquad[\mathscr{I}_{d\bar{z}^j},\mathscr{E}_{dz^k}\mathscr{E}_{d\bar{z}^l}]=g^{k\bar{j}}\mathscr{E}_{d\bar{z}^l}\]
\[[\mathscr{I}_{dz^i},L]=-\sqrt{-1}\mathscr{E}_{dz^i}\qquad[\mathscr{I}_{d\bar{z}^j},L]=\sqrt{-1}\mathscr{E}_{d\bar{z}^j}\]
\[[\mathscr{I}_{\eta},\Lambda]=0\]
\[[\mathscr{E}_{dz^i},\mathscr{I}_{dz^k}\mathscr{I}_{d\bar{z}^l}]=-g^{i\bar{l}}\mathscr{I}_{dz^k}\qquad[\mathscr{E}_{d\bar{z}^j},\mathscr{I}_{dz^k}\mathscr{I}_{d\bar{z}^l}]=g^{k\bar{j}}\mathscr{I}_{d\bar{z}^l}\]
\[[\mathscr{E}_{dz^i},\Lambda]=\sqrt{-1}\mathscr{I}_{dz^i}\qquad[\mathscr{E}_{d\bar{z}^j},\Lambda]=-\sqrt{-1}\mathscr{I}_{d\bar{z}^j}\]
\[[\mathscr{E}_{dz^i}\mathscr{E}_{d\bar{z}^j},\mathscr{I}_{dz^k}\mathscr{I}_{d\bar{z}^l}]=g^{k\bar{j}}\mathscr{E}_{dz^i}\mathscr{I}_{d\bar{z}^l}-g^{i\bar{l}}\mathscr{I}_{dz^k}\mathscr{E}_{d\bar{z}^j}\]
\[[\sqrt{-1}R^E_{i\bar{j}}\mathscr{E}_{dz^i}\mathscr{E}_{d\bar{z}^j},\Lambda]=-R_{i\bar{j}}^E\mathscr{I}_{dz^i}\mathscr{E}_{d\bar{z}^j}+R^E_{i\bar{j}}\mathscr{E}_{dz^i}\mathscr{I}_{d\bar{z}^j}=-\tr R^E+R^E_{i\bar{j}}\mathscr{E}_{d\bar{z}^j}\mathscr{I}_{dz^i}+R^E_{i\bar{j}}\mathscr{E}_{dz^i}\mathscr{I}_{d\bar{z}^j}\]
\[\tilde{\nabla}_{\partial_i}\partial_k=\partial_ig_{k\bar{q}}g^{p\bar{q}}\partial_p\qquad\tilde{\nabla}_{\partial_{\bar{j}}}\partial_{\bar{l}}=\partial_{\bar{j}}g_{r\bar{l}}g^{r\bar{s}}\partial_{\bar{s}}\]
\[\tilde{\nabla}_{\partial_i}dz^k=-\partial_ig_{p\bar{q}}g^{k\bar{q}}dz^p\qquad \tilde{\nabla}_{\partial_{\bar{j}}}d\bar{z}^l=-\partial_{\bar{j}}g_{r\bar{s}}g^{r\bar{l}}d\bar{z}^s\]
\[\mathscr{I}_{g^{p\bar{q}}\partial_p}=\mathscr{I}_{d\bar{z}^q}\qquad \mathscr{I}_{g^{p\bar{q}}\partial_{\bar{q}}}=\mathscr{I}_{dz^p}\]
\[\tilde{\nabla}_{\partial_i}dz^I=-\partial_ig_{p\bar{q}}g^{r\bar{q}}\delta^I_{rK}\,dz^p\wedge dz^K=-\partial_ig_{p\bar{q}}\mathscr{E}_{dz^p}\mathscr{I}_{d\bar{z}^q}(dz^I)=-g^{I\bar{Q}}\partial_ig_{P\bar{Q}}\,dz^P\]
\[\tilde{\nabla}_{\partial_{\bar{j}}}d\bar{z}^J=-\partial_{\bar{j}}g_{p\bar{q}}g^{p\bar{s}}\delta^{\bar{J}}_{\bar{s}\bar{L}}\,d\bar{z}^q\wedge d\bar{z}^L=-\partial_{\bar{j}}g_{p\bar{q}}\mathscr{E}_{d\bar{z}^q}\mathscr{I}_{dz^p}(d\bar{z}^J)=-g^{P\bar{J}}\partial_{\bar{j}}g_{P\bar{Q}}\,d\bar{z}^{Q}\]
Define the torsion tensor:
\[T(X,Y)=\nabla_XY-\nabla_YX-[X,Y]\]
\[T^{1,0}=\nabla_{\partial_i}\partial_k\otimes dz^i\wedge dz^k\qquad T^{0,1}=\nabla_{\partial_{\bar{j}}}\partial_{\bar{l}}\otimes d\bar{z}^j\wedge d\bar{z}^l\]
\[T_{ip\bar{q}}=\partial_ig_{p\bar{q}}-\partial_pg_{i\bar{q}}\qquad T_{\bar{j}\bar{s}r}=\partial_{\bar{j}}g_{r\bar{s}}-\partial_{\bar{s}}g_{r\bar{j}}\]
\[\mathscr{I}_{T^{1,0}}=\partial_{i}g_{p\bar{q}}\mathscr{E}_{dz^i}\mathscr{E}_{dz^p}\mathscr{I}_{d\bar{z}^q}=\frac{1}{2}T_{ip\bar{q}}\mathscr{E}_{dz^i}\mathscr{E}_{dz^p}\mathscr{I}_{d\bar{z}^q}\qquad \mathscr{I}_{T^{0,1}}=\partial_{\bar{j}}g_{r\bar{s}}\mathscr{E}_{d\bar{z}^j}\mathscr{E}_{d\bar{z}^s}\mathscr{I}_{dz^r}=\frac{1}{2}T_{\bar{j}\bar{s}r}\mathscr{E}_{d\bar{z}^j}\mathscr{E}_{d\bar{z}^s}\mathscr{I}_{dz^r}\]
For a holomorphic Hermitian vector bundle $E$ over $M$, the antilinear Hodge star operator $\bar{*}_E$ is defined by
\[\bar{*}_E:\Omega^{p,q}(E)\to\Omega^{d-p,d-q}(E^*)\]
\[\langle\psi(z),\eta(z)\rangle d\mu_M=\psi(z)\wedge\bar{*}_{E}\eta(z)\]
Then after the natural identification $E^{**}\cong E$,
\[\bar{*}_{E^*}\bar{*}_E\eta=(-1)^{p+q}\eta\]
Denote by $\tilde{\nabla}$ the total covariant derivative induced from each covariant derivative of the bundles.
\[[\mathscr{E}_{dz^i},\tilde{\nabla}_{\partial_{\bar{j}}}]=0=[\mathscr{E}_{d\bar{z}^j},\tilde{\nabla}_{\partial_i}]\]
\[[\mathscr{E}_{dz^i},\tilde{\nabla}_{\partial_k}]=-\mathscr{E}_{\tilde{\nabla}_{\partial_k}dz^i}=g^{i\bar{q}}\partial_kg_{p\bar{q}}\mathscr{E}_{dz^p}\qquad[\mathscr{E}_{d\bar{z}^j},\tilde{\nabla}_{\partial_{\bar{l}}}]=-\mathscr{E}_{\tilde{\nabla}_{\partial_{\bar{l}}}d\bar{z}^j}=g^{p\bar{j}}\partial_{\bar{l}}g_{p\bar{q}}\mathscr{E}_{d\bar{z}^q}\]
\[[\mathscr{I}_{dz^i},\tilde{\nabla}_{\partial_{\bar{j}}}]=0=[\mathscr{I}_{d\bar{z}^j},\tilde{\nabla}_{\partial_{i}}]\]
\[[\mathscr{I}_{dz^i},\tilde{\nabla}_{\partial_k}]=g^{i\bar{q}}\partial_kg_{p\bar{q}}\mathscr{I}_{dz^p}\qquad[\mathscr{I}_{d\bar{z}^j},\tilde{\nabla}_{\partial_{\bar{l}}}]=g^{p\bar{j}}\partial_{\bar{l}}g_{p\bar{q}}\mathscr{I}_{d\bar{z}^q}\]
\[[\tilde{\nabla}_X,L]=0=[\tilde{\nabla}_X,\Lambda]\]
\[\bar{*}\tilde{\nabla}_{X}\bar{*}=(-1)^{p+q}\tilde{\nabla}_{\bar{X}}\qquad \bar{*}\mathscr{E}_{\eta}\bar{*}=\mathscr{I}_{\bar{\eta}}\qquad \bar{*}\mathscr{I}_{\eta}\bar{*}=-\mathscr{E}_{\bar{\eta}}\]
\[-\bar{*}\mathscr{I}_{T^{1,0}}\bar{*}=-\partial_{\bar{j}}g_{r\bar{s}}\mathscr{I}_{d\bar{z}^j}\mathscr{I}_{d\bar{z}^s}\mathscr{E}_{dz^r}=-\frac{1}{2}T_{\bar{j}\bar{s}r}\mathscr{I}_{d\bar{z}^j}\mathscr{I}_{d\bar{z}^s}\mathscr{E}_{dz^r}=-T_{\bar{j}\bar{s}r}g^{r\bar{s}}\mathscr{I}_{d\bar{z}^j}+\frac{1}{2}T_{\bar{j}\bar{s}r}\mathscr{E}_{dz^r}\mathscr{I}_{d\bar{z}^s}\mathscr{I}_{d\bar{z}^j}\]
\[-\bar{*}\mathscr{I}_{T^{0,1}}\bar{*}=-\partial_ig_{p\bar{q}}\mathscr{I}_{dz^i}\mathscr{I}_{dz^p}\mathscr{E}_{d\bar{z}^q}=-\frac{1}
{2}T_{ip\bar{q}}\mathscr{I}_{dz^i}\mathscr{I}_{dz^p}\mathscr{E}_{d\bar{z}^q}=-T_{ip\bar{q}}g^{p\bar{q}}\mathscr{I}_{dz^i}+\frac{1}{2}T_{ip\bar{q}}\mathscr{E}_{d\bar{z}^q}\mathscr{I}_{dz^p}\mathscr{I}_{dz^i}\]
Define the Hermitian torsion:
\[\mathcal{T}=[\Lambda,\partial\Theta]=\partial_{i}g_{p\bar{q}}\mathscr{I}_{dz^i}\mathscr{E}_{dz^p}\mathscr{E}_{d\bar{z}^q}+\partial_ig_{p\bar{q}}\mathscr{E}_{dz^i}(\mathscr{I}_{dz^p}\mathscr{E}_{d\bar{z}^q}-\mathscr{E}_{dz^p}\mathscr{I}_{d\bar{z}^q})=T_{ip\bar{q}}\mathscr{I}_{dz^i}\mathscr{E}_{dz^p}\mathscr{E}_{d\bar{z}^q}-\frac{1}{2}T_{ip\bar{q}}\mathscr{E}_{dz^i}\mathscr{E}_{dz^p}\mathscr{I}_{d\bar{z}^q}\]
\[\mathcal{T}^*=T_{\bar{j}\bar{s}r}\mathscr{I}_{dz^r}\mathscr{I}_{d\bar{z}^s}\mathscr{E}_{d\bar{z}^j}-\frac{1}{2}T_{\bar{j}\bar{s}r}\mathscr{E}_{dz^r}\mathscr{I}_{d\bar{z}^s}\mathscr{I}_{d\bar{z}^j}\]
\[\bar{\mathcal{T}}=[\Lambda,\bar{\partial}\Theta]=-\partial_{\bar{j}}g_{r\bar{s}}\mathscr{I}_{d\bar{z}^j}\mathscr{E}_{dz^r}\mathscr{E}_{d\bar{z}^s}+\partial_{\bar{j}}g_{r\bar{s}}\mathscr{E}_{d\bar{z}^j}(\mathscr{I}_{d\bar{z}^s}\mathscr{E}_{dz^r}-\mathscr{E}_{d\bar{z}^s}\mathscr{I}_{dz^r})=T_{\bar{j}\bar{s}r}\mathscr{I}_{d\bar{z}^j}\mathscr{E}_{d\bar{z}^s}\mathscr{E}_{dz^r}-\frac{1}{2}T_{\bar{j}\bar{s}r}\mathscr{E}_{d\bar{z}^j}\mathscr{E}_{d\bar{z}^s}\mathscr{I}_{dz^r}\]
\[\bar{\mathcal{T}}^*=T_{ip\bar{q}}\mathscr{I}_{d\bar{z}^q}\mathscr{I}_{dz^p}\mathscr{E}_{dz^i}-\frac{1}{2}T_{ip\bar{q}}\mathscr{E}_{d\bar{z}^q}\mathscr{I}_{dz^p}\mathscr{I}_{dz^i}\]
We can write the operators $\nabla^{1,0},\bar{\partial}$ and their adjoint in terms of covariant derivative and torsions:
\[(\nabla^E)^{1,0}=\mathscr{E}_{dz^i}\tilde{\nabla}_{\partial_i}+\mathscr{I}_{T^{1,0}}=\mathscr{E}_{dz^i}\tilde{\nabla}_{\partial_i}+\frac{1}{2}T_{ip\bar{q}}\mathscr{E}_{dz^i}\mathscr{E}_{dz^p}\mathscr{I}_{d\bar{z}^q}\]
\[(\nabla^E)^{1,0*}=-\mathscr{I}_{d\bar{z}^j}\tilde{\nabla}_{\partial_{\bar{j}}}-\frac{1}{2}T_{\bar{j}\bar{s}r}\mathscr{I}_{d\bar{z}^j}\mathscr{I}_{d\bar{z}^s}\mathscr{E}_{dz^r}\]
\[\bar{\partial}^E=\mathscr{E}_{d\bar{z}^j}\tilde{\nabla}_{\partial_{\bar{j}}}+\mathscr{I}_{T^{0,1}}=\mathscr{E}_{d\bar{z}^j}\tilde{\nabla}_{\partial_{\bar{j}}}+\frac{1}{2}T_{\bar{j}\bar{s}r}\mathscr{E}_{d\bar{z}^j}\mathscr{E}_{d\bar{z}^s}\mathscr{I}_{dz^r}\]
\[\bar{\partial}^{E,*}=-\mathscr{I}_{dz^i}\tilde{\nabla}_{\partial_i}-\frac{1}
{2}T_{ip\bar{q}}\mathscr{I}_{dz^i}\mathscr{I}_{dz^p}\mathscr{E}_{d\bar{z}^q}\]

We have the following generalized Kähler identities:
\[[\bar{\partial}^{E,*},L]=\sqrt{-1}((\nabla^E)^{1,0}+\mathcal{T})=\sqrt{-1}\mathscr{E}_{dz^i}\tilde{\nabla}_{i}+\sqrt{-1}T_{ip\bar{q}}\mathscr{I}_{dz^i}\mathscr{E}_{dz^p}\mathscr{E}_{d\bar{z}^q}\]
\[[(\nabla^E)^{1,0*},L]=-\sqrt{-1}(\bar{\partial}^E+\bar{\mathcal{T}})=-\sqrt{-1}\mathscr{E}_{d\bar{z}^j}\tilde{\nabla}_{\bar{j}}-\sqrt{-1}T_{\bar{j}\bar{s}r}\mathscr{I}_{d\bar{z}^j}\mathscr{E}_{d\bar{z}^s}\mathscr{E}_{dz^r}\]
\[[\Lambda,\bar{\partial}^E]=-\sqrt{-1}((\nabla^E)^{1,0*}+\mathcal{T}^*)=\sqrt{-1}\mathscr{I}_{d\bar{z}^j}\tilde{\nabla}_{\bar{j}}+\sqrt{-1}T_{\bar{j}\bar{s}r}\mathscr{E}_{d\bar{z}^j}\mathscr{I}_{d\bar{z}^s}\mathscr{I}_{dz^r}\]
\[[\Lambda,(\nabla^E)^{1,0}]=\sqrt{-1}(\bar{\partial}^{E,*}+\bar{\mathcal{T}}^*)=-\sqrt{-1}\mathscr{I}_{dz^i}\tilde{\nabla}_i-\sqrt{-1}T_{ip\bar{q}}\mathscr{E}_{dz^i}\mathscr{I}_{dz^p}\mathscr{I}_{d\bar{z}^q}\]
The curvature of $E$:
\[R^E(X,Y)s=\nabla_X\nabla_Ys-\nabla_Y\nabla_Xs-\nabla_{[X,Y]}s\]
Let $R^E=\mathscr{E}_{dz^i}\mathscr{E}_{d\bar{z}^j}R^E_{i\bar{j}}$, then we can write
\[R^E=(\nabla^E)^2=[(\nabla^E)^{1,0},\bar{\partial}^E]\]
Define
\[\Box=[\bar{\partial},\bar{\partial}^*]\qquad \bar{\Box}=[\nabla^{1,0},\nabla^{1,0*}]\]
By Kähler identities we can get the well-known Bochner–Kodaira–Nakano identity:
\[\Box=\bar{\Box}+[\sqrt{-1}R^E,\Lambda]+[(\nabla^E)^{1,0},\mathcal{T}^*]-[\bar{\partial}^E,\bar{\mathcal{T}}^*]\]

Now for $s\in\Omega^{\bullet,\bullet}_c(M,E)$, by conisdering $\langle\langle\Box s,s\rangle\rangle$ and doing integration by part we have
\begin{equation*}
\begin{split}
||\bar{\partial}s||^2+||\bar{\partial}^*s||^2
&=||\nabla^{1,0}s||^2+||\nabla^{1,0*}s||^2+\langle\langle[\sqrt{-1}R^{E},\Lambda]s,s\rangle\rangle\\
&+\langle\langle\nabla^{1,0}s,\mathcal{T}s\rangle\rangle+\langle\langle\mathcal{T}^*s,\nabla^{1,0*}s\rangle\rangle+\langle\langle\bar{\partial}s,\bar{\mathcal{T}}s\rangle\rangle+\langle\langle\bar{\mathcal{T}}^*s,\bar{\partial}^*s\rangle\rangle
\end{split}
\end{equation*}
Then by Cauchy inequaltiy, we get the Bochner–Kodaira–Nakano inequality
\begin{equation}\label{BKN inequality}
\begin{split}
\frac{3}{2}(||\bar{\partial}s||^2+||\bar{\partial}^*s||^2)
&\ge\frac{1}{2}(||\nabla^{1,0}s||^2+||\nabla^{1,0*}s||^2)+\langle\langle[\sqrt{-1}R^{E},\Lambda]s,s\rangle\rangle\\
&-||\mathcal{T}s||^2-||\mathcal{T}^*s||^2-||\bar{\mathcal{T}}s||^2-||\bar{\mathcal{T}}^*s||^2
\end{split}
\end{equation}
Notice that if $M$ is Kähler, then $\mathcal{T}=0$ and we have a slightly better constant.

Next we derive the formula for manifold with smooth boundary. For any $z\in\partial M$, we can locally take a smooth real-valued function $\rho$ defined on a coordinate neighborhood such that the interior of $M$ is described by $\rho<0$ and $|d\rho|=1$ on $\partial M$. Let $\Omega^{\bullet,\bullet}_c(\bar{M},E)$ be the forms smooth up to boundary and compactly supported in $\bar{M}=M\cup\partial M$.
We have the basic formulas for $s_1,s_2\in\Omega^{\bullet,\bullet}_c(\bar{M},E)$:
\[\langle\langle\bar{\partial}s_1,s_2\rangle\rangle-\langle\langle s_1,\bar{\partial}^*s_2\rangle\rangle=\int_{\partial M}\langle\bar{\partial}\rho\wedge s_1,s_2\rangle\]
\[\langle\langle\nabla^{1,0}s_1,s_2\rangle\rangle-\langle\langle s_1,\nabla^{1,0*}s_2\rangle\rangle=\int_{\partial M}\langle\partial\rho\wedge s_1,s_2\rangle\]

Now for $s\in\Omega^{\bullet,\bullet}_c(\bar{M},E)$ we have
\begin{equation*}
\begin{split}
||\bar{\partial}s||^2+||\bar{\partial}^*s||^2
&=||\nabla^{1,0}s||^2+||\nabla^{1,0*}s||^2+\langle\langle[\sqrt{-1}R^{E},\Lambda]s,s\rangle\rangle\\
&+\langle\langle\nabla^{1,0}s,\mathcal{T}s\rangle\rangle+\langle\langle\mathcal{T}^*s,\nabla^{1,0*}s\rangle\rangle+\langle\langle\bar{\partial}s,\bar{\mathcal{T}}s\rangle\rangle+\langle\langle\bar{\mathcal{T}}^*s,\bar{\partial}^*s\rangle\rangle\\
&+\int_{\partial M}\langle\bar{\partial}s,\bar{\partial}\rho\wedge s\rangle+\int_{\partial M}\langle\partial\rho\wedge(\nabla^{1,0*}+\mathcal{T}^*)s,s\rangle\\
&-\int_{\partial M}\langle\bar{\partial}\rho\wedge(\bar{\partial}^*+\bar{\mathcal{T}}^*)s,\Psi^F\wedge s\rangle-\int_{\partial M}\langle\nabla^{1,0}s,\partial\rho\wedge s\rangle\\
\end{split}
\end{equation*}
Again we have the inequality
\begin{equation}\label{BKN inequality with boundary terms}
\begin{split}
\frac{3}{2}(||\bar{\partial}s||^2+||\bar{\partial}^*s||^2)
&\ge\frac{1}{2}(||\nabla^{1,0}s||^2+||\nabla^{1,0*}s||^2)+\langle\langle[\sqrt{-1}R^{E},\Lambda]s,s\rangle\rangle\\
&-||\mathcal{T}s||^2-||\mathcal{T}^*s||^2-||\bar{\mathcal{T}}s||^2-||\bar{\mathcal{T}}^*s||^2\\
&+\int_{\partial M}\langle\bar{\partial}s,\bar{\partial}\rho\wedge s\rangle+\int_{\partial M}\langle \partial\rho\wedge(\nabla^{1,0*}+\mathcal{T}^*)s,s\rangle\\
&-\int_{\partial M}\langle \bar{\partial}
\rho\wedge(\bar{\partial}^*+\bar{\mathcal{T}}^*)s,\Psi^F\wedge s\rangle-\int_{\partial M}\langle\nabla^{1,0}s,\partial\rho\wedge s\rangle\\
\end{split}
\end{equation}

Define the Dirac-Dolbeault operator
\[D=\sqrt{2}(\bar{\partial}+\bar{\partial}^*)\]
Then $\frac{1}{2}D^2=\Box$ and for $s\in\Omega^k(M,E)$,
\[\frac{1}{2}|Ds|^2=|\bar{\partial}s|^2+|\bar{\partial}^*s|^2\]

Now for $\eta_1\in\Omega^{k_1}(M,E_1),\eta_2\in\Omega^{k_2}(M,E_2)$,
\[\bar{\partial}^{E_1\otimes E_2}(\eta_1\wedge\eta_2)=\bar{\partial}^{E_1}\eta_1\wedge\eta_2+(-1)^{k_1}\eta_1\wedge\bar{\partial}^{E_2}\eta_2\]
\begin{equation*}
\begin{split}
\bar{\partial}^{E_1\otimes E_2,*}(\eta_1\wedge\eta_2)
&=\bar{\partial}^{E_1,*}\eta_1\wedge\eta_2+(-1)^{k_1}\eta_1\wedge\bar{\partial}^{E_2,*}\eta_2\\
&-\mathscr{I}_{dz^i}\eta_1\wedge(\tilde{\nabla}_i-\frac{1}{2}T_{ip\bar{q}}\mathscr{E}_{d\bar{z}^q}\mathscr{I}_{dz^p})\eta_2-(-1)^{k_1}(\tilde{\nabla}_i-\frac{1}{2}T_{ip\bar{q}}\mathscr{E}_{d\bar{z}^q}\mathscr{I}_{dz^p})\eta_1\wedge\mathscr{I}_{dz^i}\eta_2
\end{split}
\end{equation*}
For the special case $\eta_1\in\Omega^{p,0}(M,E_1),\eta_2\in\Omega^{0,q}(M,E_2)$,
\[\bar{\partial}(\eta_1\wedge\eta_2)=(-1)^p(\tilde{\nabla}_{\bar{j}}\eta_1\wedge\mathscr{E}_{d\bar{z}^j}\eta_2+\eta_1\wedge\bar{\partial}\eta_2)\]
\[\bar{\partial}^*(\eta_1\wedge\eta_2)=(-1)^p(-\tilde{\nabla}_i\eta_1\wedge\mathscr{I}_{dz^i}\eta_2+\eta_1\wedge\bar{\partial}^*\eta_2)\]
Thus
\begin{equation}\label{Leibniz rule}
\frac{1}{2}|D(\eta_1\wedge\eta_2)|^2=|\tilde{\nabla}_{\bar{j}}\eta_1\wedge\mathscr{E}_{d\bar{z}^j}\eta_2+\eta_1\wedge\bar{\partial}\eta_2|^2+|-\tilde{\nabla}_i\eta_1\wedge\mathscr{I}_{dz^i}\eta_2+\eta_1\wedge\bar{\partial}^*\eta_2|^2
\end{equation}

\section{Linear algebraic constructions}\label{section:linear algebra}

In this section, we construct the section we are going to twist in an abstract setting, derive several properties, and compute a quantity that corresponds to the error term.

\begin{definition}
Let $E$ be a Hermitian vector bundle over a manifold $M$, and $F\subset E$ be a smooth $1$-dimensional subbundle. Locally take a nonvanishing section $\psi$ and define
\[\Psi^F=|\psi|^{-2}\psi\otimes\psi^*\in C^{\infty}(M,E\otimes E^*)\]
where $\psi^*$ is the metric dual. This definition is independent of the choice of $\psi$ and can be defined globally.
\end{definition}

\begin{proposition}\label{properties of characteristic section}
Let $\nabla^E$ be a metric connection on $E$, naturally extend to $E\otimes E^*$, we have
\begin{enumerate}[(a)]
    \item $\langle\nabla^E_X\Psi^F,\Psi^F\rangle=0$. As a result, also $\langle\nabla^E_{\bar{Y}}\nabla^E_X\Psi^F,\Psi^F\rangle=-\langle\nabla^E_X\Psi^F,\nabla^E_{Y}\Psi^F\rangle$.
    \item $\langle R^{E\otimes E^*}\Psi^F,\Psi^F\rangle=0$.
    \item $\langle\nabla^E_X\Psi^F,\nabla^E_Y\Psi^F\rangle=\langle\nabla^E_{\bar{Y}}\Psi^F,\nabla^E_{\bar{X}}\Psi^F\rangle$
\end{enumerate}
\begin{proof}
\begin{enumerate}[(a)]
    \item \[\langle\nabla^E_X\Psi^F,\Psi^F\rangle=|\psi|^{-2}d_X(|\psi|^{-2})\langle\psi\otimes\psi^*,\psi\otimes\psi^*\rangle+|\psi|^{-4}\langle\nabla^E_X\psi\otimes\psi^*+\psi\otimes\nabla^E_X\psi^*,\psi\otimes\psi^*\rangle\]
    Notice that
    \[\langle\nabla^E_X\psi^*,\psi^*\rangle=\langle(\nabla^E_{\bar{X}}\psi)^*,\psi^*\rangle=\langle\psi,\nabla^E_{\bar{X}}\psi\rangle\]
    Thus,
    \[\langle\nabla^E_X\Psi^F,\Psi^F\rangle=|\psi|^2d_X(|\psi|^{-2})+|\psi|^{-2}(\langle\nabla^E_X\psi,\psi\rangle+\langle\psi,\nabla^E_{\bar{X}}\psi\rangle)=d_X(|\psi|^2\cdot|\psi|^{-2})=0\]
    \item \[\langle R^{E\otimes E^*}(X,Y)\Psi^F,\Psi^F\rangle=|\psi|^{-2}(\langle R^E(X,Y)\psi,\psi\rangle+\langle R^{E^*}(X,Y)\psi^*,\psi^*\rangle)\]
    Notice that
    \[\langle R^{E^*}(X,Y)\psi^*,\psi^*\rangle=\langle (R^E(\bar{X},\bar{Y})\psi)^*,\psi^*\rangle=\langle\psi,R^E(\bar{X},\bar{Y})\psi\rangle\]
    Apply the operator $XY-YX-[X,Y]=0$ to $\langle\psi,\psi\rangle$ and apply metric compatibility we get
    \[\langle R^E(X,Y)\psi,\psi\rangle+\langle\psi,R^E(\bar{X},\bar{Y})\psi\rangle=0\]
    \item By (a),
    \[\langle\nabla^E_X\Psi^F,\nabla^E_{Y}\Psi^F\rangle=-\langle\nabla^E_{\bar{Y}}\nabla^E_X\Psi^F,\Psi^F\rangle\]
    \[\langle\nabla^E_{\bar{Y}}\Psi^F,\nabla^E_{\bar{X}}\Psi^F\rangle=-\langle\nabla^E_X\nabla^E_{\bar{Y}}\Psi^F,\Psi^F\rangle\]
    and by (b) RHS are equal.
\end{enumerate}
\end{proof}
\end{proposition}

\begin{proposition}
Assume there is a bundle map $A$ on $E$ that is fibre-wise self-adjoint, and its maximal eigenvalue has multiplicity one at every point, then the corresdponding eigenspace $F$ forms a $1$-dimensional smooth subbundle of $E$.
\begin{proof}
This is a local problem, we only need to construct a local nonvanishing smooth section near a fixed point $z_0\in M$. Take a local frame $\{e^\gamma\}$ of $E$ (In the later theorem, we will consider $E=\bigwedge\nolimits^{p,0}$, so write upper indices 
just for temporary consistency). We pick such frame so that $\{e^\gamma(z_0)\}$ is a orthonormal eigenbasis with $Ae^\alpha=\lambda_1(z_0)e^\alpha$ obtaining the maximal eigenvalue $\lambda_1(z)$ at $z_0$. Represent $A=(a_{\gamma}^{\beta})$ in the basis $\{e^\gamma\}$. We place $a_\alpha^\alpha$ to the upper left corner and write the matrix $A$ in the block form
\[A=\matrixtwo{a_1^1}{a_1^2}{a_2^1}{a_2^2}\]
where 
\[a_1^1=a_\alpha^\alpha,\ a_1^2=\left( \begin{array}{ccc}
- & a_\alpha^\beta &-
\end{array} \right)_{\beta\neq \alpha},\ a_2^1=\left( \begin{array}{c}
| \\
a_\beta^\alpha\\
|
\end{array} \right)_{\beta\neq \alpha},\ a_2^2=(a_\gamma^\delta)_{\gamma,\delta\neq \alpha}\]
We can consider solving $v(z)$ in the equation
\[\matrixtwo{a_1^1}{a_1^2}{a_2^1}{a_2^2}\vectortwo{1}{v(z)}=\lambda_1(z)\vectortwo{1}{v(z)}\]
with the condition $v(z_0)=0$. By continuity, near $z_0$ every nonzero maximal eigenvector should have the coefficient of $e^{\alpha}$ nonzero, by the multiplicity one assumption and normalizing we see that there exists a unique $v(z)$ satisfying this equation for $z$ near $z_0$. Now
\[(\lambda_1(z)-a_2^2)v(z)=a_2^1\]
By uniqueness of $v(z)$, we also know that $\lambda_1(z)-a_2^2$ is invertible, and we can write
\[v(z)=(\lambda_1(z)-a_2^2)^{-1}a_2^1\]
From this expression we know $v(z)$ is smooth, and this gives a smooth nonvanishing local section of $F$.
\end{proof}
\end{proposition}

\begin{thm}
Let $(M,\Theta)$ be a Hermitian manifold and assume there is a bundle map $A$ on $\bigwedge\nolimits^{p,0}$ that is fibrewise self-adjoint. Suppose that the maximal eigenvalue of $A$ has multiplicity one at the point $z_0\in M$. Represent $A=(a_K^L)_{K,L},\ K,L\in\binom{[d]}{p}$ with respect to the basis $\{dz^K\}$. Further assume that there exists a coordinate neighborhood $(U,z)$ of $z_0$ such that
\[g_{i\bar{j}}(z_0)=\delta_{ij},\ a_K^L(z_0)=\delta_K^L\lambda_K\]
Let $I$ be the index that $\lambda_I=\lambda_1(z_0)$ is maximal, and let $F$ be the smooth $1$-dimensional subbundle correspond to maximal eigenvalue $\lambda_1(z)$ near $z_0$. Then, we have the computation
\[|\bar{\partial}\Psi^F|^2|_{z=z_0}=\sum_{J\neq I}(\lambda_I-\lambda_J)^{-2}(|\bar{\partial}a_I^J|^2+|\bar{\partial}a_J^I|^2)\]
\begin{proof}
We place the index $I$ to the upper left corner and write the matrix $A$ in the block form
\[A=\matrixtwo{a_1^1}{a_1^2}{a_2^1}{a_2^2}\]
where 
\[a_1^1=a_I^I,\ a_1^2=\left( \begin{array}{ccc}
- & a_I^J &-
\end{array} \right)_{J\neq I},\ a_2^1=\left( \begin{array}{c}
| \\
a_J^I\\
|
\end{array} \right)_{J\neq I},\ a_2^2=(a_K^L)_{K,L\neq I}\]
We can also write $g=(g_{K\bar{L}})$ into the same block form.
As in the proof of previous proposition, we can take a local section $\psi$ of $F$ to be
\[\psi=dz^I+\left( \begin{array}{ccc}
- & dz^K &-
\end{array} \right)_{K\neq I}(\lambda_1-a_2^2)^{-1}a_2^1\]
\[\psi^*=\langle dz^J,\psi\rangle\partial_J=(g^{J\bar{I}}+\left( \begin{array}{ccc}
- & g^{J\bar{K}} &-
\end{array} \right)_{K\neq I}\overline{(\lambda_1-a_2^2)^{-1}a_2^1})\partial_J\]
\begin{equation*}
\begin{split}
|\psi|^2
&=g^{I\bar{I}}+\left( \begin{array}{ccc}
-&g^{I\bar{K}}&-
\end{array} \right)_{K\neq I}\overline{(\lambda_1-a_2^2)^{-1}a_2^1}+\left( \begin{array}{ccc}
-&g^{K\bar{I}}&-
\end{array} \right)_{K\neq I}(\lambda_1-a_2^2)^{-1}a_2^1\\
&+(a_2^1)^t(\lambda_1-a_2^2)^{-t}(g^{K\bar{L}})_{K,L\neq I}\overline{(\lambda_1-a_2^2)^{-1}a_2^1}
\end{split}
\end{equation*}

\begin{equation*}
\begin{split}
\psi\otimes\psi^*
&=(g^{J\bar{I}}+\left( \begin{array}{ccc}
- & g^{J\bar{K}} &-
\end{array} \right)_{K\neq I}\overline{(\lambda_1-a_2^2)^{-1}a_2^1})\,dz^I\otimes\partial_J\\
&+\left( \begin{array}{ccc}
- & dz^J\otimes\partial_K &-
\end{array} \right)_{J\neq I}(\lambda_1-a_2^2)^{-1}a_2^1(g^{K\bar{I}}+\left( \begin{array}{ccc}
- & g^{K\bar{L}} &-
\end{array} \right)_{L\neq I}\overline{(\lambda_1-a_2^2)^{-1}a_2^1})
\end{split}
\end{equation*}
Notice that $a_1^2(z_0)=0$, $a_2^1(z_0)=0$ and $g^{K\bar{L}}(z_0)=\delta^{KL}$, we can compute
\begin{equation*}
\begin{split}
\bar{\partial}\Psi^F|_{z=z_0}
&=\bar{\partial}|\psi|^{-2}|_{z=z_0}dz^I\otimes\partial_I+\bar{\partial}(\psi\otimes\psi^*)|_{z=z_0}\\
&=\bar{\partial}g_{I\bar{I}}|_{z=z_0}dz^I\otimes\partial_I+\sum_J(-\bar{\partial}g_{I\bar{J}}+\left( \begin{array}{ccc}
- & \delta^{JK} &-
\end{array} \right)_{K\neq I}\overline{(\lambda_{1}-a^2_2)^{-1}}\bar{\partial}\overline{a^1_2})\,dz^I\otimes\partial_J\\
&+\left( \begin{array}{ccc}
-&dz^J\otimes\partial_I&-
\end{array} \right)_{J\neq I}(\lambda_1-a_2^2)^{-1}\bar{\partial}a_2^1\\
&=\sum_{J\neq I}(-\bar{\partial}g_{I\bar{J}}+\left( \begin{array}{ccc}
- & \delta^{JK} &-
\end{array} \right)_{K\neq I}\overline{(\lambda_{1}-a^2_2)^{-1}}\bar{\partial}\overline{a^1_2})\,dz^I\otimes\partial_J\\
&+\left( \begin{array}{ccc}
-&dz^J\otimes\partial_I&-
\end{array} \right)_{J\neq I}(\lambda_1-a_2^2)^{-1}\bar{\partial}a_2^1\\
&=\sum_{J\neq I}(-\bar{\partial}g_{I\bar{J}}+\left( \begin{array}{ccc}
- & dz^I\otimes\partial_J &-
\end{array} \right)_{J\neq I}\overline{(\lambda_{1}-a^2_2)^{-1}}\bar{\partial}\overline{a^1_2})\\
&+\left( \begin{array}{ccc}
-&dz^J\otimes\partial_I&-
\end{array} \right)_{J\neq I}(\lambda_1-a_2^2)^{-1}\bar{\partial}a_2^1\\
\end{split}
\end{equation*}
Computing the norm, we directly get
\begin{equation*}
\begin{split}
|\bar{\partial}\Psi^F|^2|_{z=z_0}
&=\sum_{J\neq I}|-\bar{\partial}g_{I\bar{J}}+(\lambda_I-\lambda_J)^{-1}\bar{\partial}\overline{a_J^I}|^2+\sum_{J\neq I}|(\lambda_I-\lambda_J)^{-1}\bar{\partial}a_J^I|^2
\end{split}
\end{equation*}

Now notice that $A$ is self-adjoint operator, which means that $\overline{a_{\alpha\bar{\beta}}}=a_{\beta\bar{\alpha}}$, so that
\[\overline{a_J^I}=\overline{a_{J\bar{\alpha}}g^{I\bar{\alpha}}}=a_{\alpha\bar{J}}g^{\alpha\bar{I}}=a_{\alpha}^\beta g_{\beta\bar{J}}g^{\alpha\bar{I}}\]
\[\bar{\partial}\overline{a_J^I}=\bar{\partial}a_{\alpha}^\beta g_{\beta\bar{J}}g^{\alpha\bar{I}}+a_{\alpha}^\beta\bar{\partial}g_{\beta\bar{J}}g^{\alpha\bar{I}}-a_{\alpha\beta}g_{\beta\bar{J}}\bar{\partial}g_{I\bar{\alpha}}\]
Thus at $z=z_0$ we have
\[\bar{\partial}\overline{a_J^I}|_{z=z_0}=\bar{\partial}a_{I}^J+\lambda_I\bar{\partial}g_{I\bar{J}}-\lambda_J\bar{\partial}g_{I\bar{J}}\]
As a result, we finish the computation
\[|\bar{\partial}\Psi^F|^2|_{z=z_0}=\sum_{J\neq I}(\lambda_I-\lambda_J)^{-2}(|\bar{\partial}a_I^J|^2+|\bar{\partial}a_J^I|^2)\]

\end{proof}
\end{thm}
    
\section{$L^2$ estimates for Dirac-Dolbeault operator}\label{section:L^2 estimate}
Denote by $\lambda^L_1(z)\ge\lambda^L_2(z)\ge\cdots\ge\lambda^L_d(z)$ the eigenvalues of $R^L(z)$ and $D_k$ the Dirac-Dolbeault operator with coefficient $E\otimes L^k$.

\subsection{Main result}

\begin{thm}[Main theorem]\label{thm:Main theorem}
Let $(M,\Theta)$ be a Hermitian manifold of dimension $d$. $E$ and $L$ are holomorphic vector bundle and line bundle with Hermitian metrics.
Assume that there is a smooth $1$-dimensional subbundle $F$ of $\bigwedge\nolimits^{p,0}$ such that for a local section $\psi(z)$ of $F$ we have the estimate
\[\langle[\sqrt{-1}R^L,\Lambda]\psi(z),\psi(z)\rangle\ge \lambda^F(z)|\psi|^2(z)\]

Assume $[\sqrt{-1}R^{E\otimes\bigwedge\nolimits_{p,0}},\Lambda](z)\ge-\lambda^E_p(z)\id^{E\otimes \bigwedge\nolimits_{p,0}}(z)$ on $\bigwedge\nolimits^{p,q}$. Then, for any $s\in\Omega^{0,q}_c(M,E\otimes L^k)$ we have the inequality
\[\frac{1}{2}||D_ks||^2\ge \int_M(\frac{2}{3}k(\lambda^F+\sum_{j=d-q+1}^{d}\lambda^L_j)-\frac{2}{3}\lambda^E_p-C_d|\partial\Theta|^2-|\bar{\partial}\Psi^F|^2)|s|^2\]
where $C_d$ is a constant only depending on dimension $d$.

If moreover that $(M,\Theta)$ is Kähler, we have
\[\frac{1}{2}||D_ks||^2\ge\int_M (k(\lambda^F+\sum_{j=d-q+1}^{d}\lambda^L_j)-\lambda^E_p-|\bar{\partial}\Psi^F|^2)|s|^2\]

\begin{proof}
Apply the Bochner-Kodaira-Nakano inequality (\ref{BKN inequality}) to $\Psi^F\wedge s\in\Omega^{p,q}_c(M,E\otimes\bigwedge\nolimits_{p,0}\otimes L^k)$:
\begin{equation}\label{eq:1}
\begin{split}
\frac{1}{2}||D_k(\Psi^F\wedge s)||^2
&\ge\frac{2}{3}\langle\langle[\sqrt{-1}R^{E\otimes\bigwedge\nolimits_{p,0}\otimes L^k},\Lambda]\Psi^F\wedge s,\Psi^F\wedge s\rangle\rangle\\
&-\frac{1}{3}(||\mathcal{T}(\Psi^F\wedge s)||^2+||\mathcal{T}^*(\Psi^F\wedge s)||^2+||\bar{\mathcal{T}}(\Psi^F\wedge s)||^2+||\bar{\mathcal{T}}^*(\Psi^F\wedge s)||^2)
\end{split}
\end{equation}
By assumption
\[\langle[\sqrt{-1}R^{E\otimes\bigwedge\nolimits_{p,0}\otimes L^k},\Lambda]\Psi^F\wedge s,\Psi^F\wedge s\rangle\ge (k(\lambda^F+\sum_{j=d-q+1}^{d}\lambda^L_j)-\lambda^E_p)|s|^2\]
Also we can find a $C_d$ with
\[\frac{1}{3}(|\mathcal{T}|^2_{\text{op}}+|\mathcal{T}^*|^2_{\text{op}}+|\bar{\mathcal{T}}|^2_{\text{op}}+|\bar{\mathcal{T}}^*|^2_{\text{op}})\le C_d|\partial\Theta|^2\]
Finally, from (\ref{Leibniz rule}):
\[\frac{1}{2}|D_k(\Psi^F\wedge s)|^2=|\tilde{\nabla}_{\bar{j}}\Psi^F\wedge\mathscr{E}_{d\bar{z}^j}s+\Psi^F\wedge\bar{\partial}s|^2+|-\tilde{\nabla}_i\Psi^F\wedge\mathscr{I}_{dz^i}s+\Psi^F\wedge\bar{\partial}^*s|^2\]
Using the properties (\ref{properties of characteristic section})
\[\langle\tilde{\nabla}_X\Psi^F,\Psi^F\rangle=0\qquad \langle\tilde{\nabla}_{i}\Psi^F,\tilde{\nabla}_{j}\Psi^F\rangle=\langle\tilde{\nabla}_{\bar{j}}\Psi^F,\tilde{\nabla}_{\bar{i}}\Psi^F\rangle\]
and the basic equation $[\mathscr{I}_{dz^i},\mathscr{E}_{d\bar{z}^j}]=g^{i\bar{j}}$, we get
\[\frac{1}{2}|D_k(\Psi^F\wedge s)|^2=\frac{1}{2}|D_ks|^2+g^{i\bar{j}}\langle\tilde{\nabla}_{\bar{j}}\Psi^F,\tilde{\nabla}_{\bar{i}}\Psi^F\rangle|s|^2=\frac{1}{2}|D_ks|^2+|\bar{\partial}\Psi^F|^2|s|^2\]
This equality holds pointwise and in later computation with boundary (\ref{L^2 estimates with boundary}) it will not contribute to boundary terms.

Insert all equations back to (\ref{eq:1}) we get the result. 

If $(M,\Theta)$ is Kähler, we have a better Bochner Kodaira inequality:
\[\frac{1}{2}||D_k(\Psi^F\wedge s)||^2
\ge\langle\langle[\sqrt{-1}R^{E\otimes\bigwedge\nolimits_{p,0}\otimes L^k},\Lambda]\Psi^F\wedge s,\Psi^F\wedge s\rangle\rangle\]
And we get the corresponding better result.

\end{proof}

\end{thm}

Now combine this with the discussion in previous section.

Assume $\lambda_p>\lambda_{p+1}$, then there exists a smooth $1$-dimensional subbundle $F_p$ of $\bigwedge\nolimits^{p,0}$ such that its fibres are the maximal eigenspace of $[\sqrt{-1}R^L,\Lambda]$ (with eigenvalue $\sum_{i=1}^p\lambda^L_i-\sum_{k=1}^d \lambda^L_k$). 

Denote $\lambda_{p,q}^L=\sum_{i=1}^p\lambda^L_i-\sum_{j=1}^{d-q}\lambda^L_{j}$, then for $s\in\Omega^{0,q}_c(M,E\otimes L^k)$, we have
\[\langle[\sqrt{-1}R^{E\otimes\bigwedge\nolimits_{p,0}\otimes L^k},\Lambda]\Psi^F\wedge s,\Psi^F\wedge s\rangle\ge (k\lambda_{p,q}-\lambda_p^E)|s|^2\]

Thus we obtain the following theorem:
\begin{thm}\label{thm:Main theorem applied to important bundle}
Let $(M,\Theta)$ be a Hermitian manifold of dimension $d$. $E$ and $L$ are holomorphic vector bundle and line bundle with Hermitian metrics. Assume $\lambda_p>\lambda_{p+1}$ and $[\sqrt{-1}R^{E\otimes\bigwedge\nolimits_{p,0}},\Lambda](z)\ge-\lambda^E_p(z)\id^{E\otimes \bigwedge\nolimits_{p,0}}(z)$ on $\bigwedge\nolimits^{p,q}$, then we have the following estimate for Dirac-Dobeault operator $D_k$:
\[\frac{1}{2}||D_ks||^2\ge \int_M\prescript{}{k}{f}_{p,q}|s|^2\]
for all $s\in\Omega^{0,q}_c(M,E\otimes L^k)$, where 
\[\prescript{}{k}{f}_{p,q}=
\begin{cases}
k\lambda_{p,q}-\lambda^E_p-|\bar{\partial}\Psi^{F_p}|^2 &\text{ if $M$ is Kähler}\\
\frac{2}{3}k\lambda_{p,q}-\frac{2}{3}\lambda^E_p-C_d|\partial\Theta|^2-|\bar{\partial}\Psi^{F_p}|^2 &\text{ if $M$ is non-Kähler}
\end{cases}\]
and $C_d$ is a constant only depending on dimension $d$.

\end{thm}

\begin{remark}
For $k>>0$, we usually need $\lambda_{p,q}\ge0$ to deduce something.

Notice that if $p+q>d$, then $\lambda_{p,q}\ge0$ implies $\lambda_1\ge\cdots\ge\lambda_{d-q+1}\ge0$. If $p+q<d$, then $\lambda_{p,q}\ge0$ implies $0\ge\lambda_{d-q}\ge\cdots\ge\lambda_d$. 
\end{remark}

\begin{remark}
This theorem has the multiplicity assumption $\lambda_p>\lambda_{p+1}$ to ensure that $F_p$ is smooth. This condition is naturally satisfied if the curvature $R^L$ is nondegenerate with signature $(d-p,p)$ everywhere.

If we do not pick the fibre to be the maximal eigensubspace of $[\sqrt{-1}R^L,\Lambda]$, we can take other smooth $1$-dimensional subbundle $F$. We can drop the multiplicity assumption if the desired region is compact and $\lambda_{p,q}>0$ there, and we will present this later. However, to deal with those points $\lambda_{p,q}=0$ (degenerate points) possibly the assumption is technically needed in this method.
\end{remark}

Next we give some bounds for the quantity $|\bar{\partial}\Psi^{F_p}|^2$ in this case.
\begin{proposition}
If $\lambda_p>\lambda_{p+1}$ and denote
\[R^L=R_{i\bar{j}}\,dz^i\wedge d\bar{z}^j\]
Let $R:T^{1,0}\to T^{1,0}$ be its type change to $(1,1)$-tensor by
\[R=R_i^j\,dz^i\otimes\partial_j=R_{i\bar{s}}g^{j\bar{s}}\,dz^i\otimes\partial_j\]
If we pick a local coordinate near $z_0\in M$ such that $R$ is unitarily diagonalized by $\{\partial_i|_{z=z_0}\}_{i=1}^d$ with eigenvalues $\lambda_1\ge\cdots\ge\lambda_d$, then
\[|\bar{\partial}\Psi^{F_p}|^2|_{z=z_0}=\sum_{\substack{i\le p\\ j\ge p+1}}(\lambda_i-\lambda_j)^{-2}(|\bar{\partial}R^i_j|^2+|\bar{\partial}R_i^j|^2)\]
Moreover, we have the estimate independent of coordinate choice:
\[|\bar{\partial}\Psi^{F_p}|^2\le(\lambda_p-\lambda_{p+1})^{-2}|\bar{\partial}R|^2\]
Also, there are two slightly better estimates:
\begin{enumerate}[(a)]
    \item \[|\bar{\partial}\Psi^{F_p}|^2\le(\lambda_p-\lambda_{p+1})^{-2}(|\bar{\partial}R|^2-\frac{1}{d}|\bar{\partial}\tr R|^2)\]
    \item For points $|R|\neq0$,
    \[|\bar{\partial}\Psi^{F_p}|^2\le(\lambda_p-\lambda_{p+1})^{-2}(|\bar{\partial}R|^2-|\bar{\partial}|R||^2)\]
\end{enumerate}
\begin{proof}
\[[\sqrt{-1}R_{i\bar{j}}\mathscr{E}_{dz^i}\mathscr{E}_{d\bar{z}^j},\Lambda]=-\tr R^L+ R_{i\bar{j}}\mathscr{E}_{d\bar{z}^j}\mathscr{I}_{dz^i}+R_{i\bar{j}}\mathscr{E}_{dz^i}\mathscr{I}_{d\bar{z}^j}\]
The part $-\tr R^L$ does not affect $\Psi^{F_p}$, and $\mathscr{I}_{dz^i}=0$ on $\bigwedge\nolimits^{p,0}$, so we only need to see how the operator $A=R_{i\bar{j}}\mathscr{E}_{dz^i}\mathscr{I}_{d\bar{z}^j}$ acts on $\bigwedge\nolimits^{p,0}$. We assume that the coordinate satisfies $g_{i\bar{j}}(z_0)=\delta_{i\bar{j}}$ and for $I=[p]$, $dz^I|_{z=z_0}$ obtain the maximal eigenvalue. So
\[a_I^J=R_i^j\delta^J_{j K}\delta^{i K}_{I}\]
Thus by linear algebra,
\begin{equation*}
\begin{split}
|\bar{\partial}\Psi^{F_p}|^2|_{z=z_0}
&=\sum_{J\neq I}(\lambda_I-\lambda_J)^{-2}(|\bar{\partial}a_I^J|^2+|\bar{\partial}a_J^I|^2)\\
&=\sum_{I\cap J=K\in\binom{[d]}{p-1}}(\lambda_I-\lambda_J)^{-2}(|\delta_{jK}^J\delta_I^{iK}\bar{\partial}R_i^j|^2+|\delta_{iK}^I\delta_J^{jK}\bar{\partial}R_j^i|^2)\\
&=\sum_{i\in I,j\not\in I}(\lambda_i-\lambda_j)^{-2}(|\bar{\partial}R_i^j|^2+|\bar{\partial}R_j^i|^2)
\end{split}
\end{equation*}

Notice that $\lambda_1\ge\cdots\ge\lambda_p>\lambda_{p+1}\ge\cdots\ge\lambda_d$ and $I=[p]$, so we have $\lambda_i-\lambda_j\ge\lambda_p-\lambda_{p+1}$ for $i\in I,j\not\in I$. Thus,
\[\sum_{i\in I,j\not\in I}(\lambda_i-\lambda_j)^{-2}(|\bar{\partial}R_i^j|^2+|\bar{\partial}R_j^i|^2)\le(\lambda_p-\lambda_{p+1})^{-2}(\sum_{i,j,k}|\partial_{\bar{k}}R_i^j|^2-\sum_{i,k}|\partial_{\bar{k}}R_i^i|^2)\]
Notice that $\sum_{i,j,k}|\partial_{\bar{k}}R_i^j|^2=|\bar{\partial}R|^2$. We can obtain a bound if we obtain a lower bound for $\sum_{i,k}|\partial_{\bar{k}}R_i^i|^2$ (a trivial lower bound is $0$).
\begin{enumerate}[(a)]
    \item Cauchy-Schwarz inequality:
    \[|\bar{\partial}\tr R|^2(z_0)=\sum_k|\sum_{i}\partial_{\bar{k}}R_i^i|^2\le d\sum_{i,k}|\partial_{\bar{k}}R_i^i|^2\]
    Thus,
    \[\sum_{i,k}|\partial_{\bar{k}}R_i^i|^2\ge\frac{1}{d}|\bar{\partial}\tr R|^2\]
    \item 
    Again by Cauchy-Schwarz:
    \[|\bar{\partial}|R^L||^2(z_0)=|R^L|^{-2}\sum_k|\sum_{i}\partial_{\bar{k}}R_i^iR_i^i|^2\le\sum_{i,k}|\partial_{\bar{k}}R_i^i|^2\]
\end{enumerate}

\end{proof}
\end{proposition}

Next, we apply the main theorem to a much more trivial case.

\begin{proposition}
For every point $z_0\in M$ with $\lambda_{p,q}(z_0)>0$, and for any $\epsilon>0$ there exists a open set $U$ containing $z_0$ so that for $s\in\Omega^{0,q}_c(U,E\otimes L^k)$, we have
\[\frac{1}{2}||D_ks||^2\ge (k(\lambda_{p,q}(z_0)-\epsilon)-C)||s||^2\]
\begin{proof}
Take a coordinate neigborhood $U$ of $z_0$ and at the point $z_0$ we unitarily diagonalize $R^L(z_0)$ and assume that $I\subset[d]$ satisfies for all $\eta\in\Omega^{0,q}_c(U,E\otimes L^k)$,
\[\langle[\sqrt{-1}R^L,\Lambda](dz^I\wedge\eta),dz^I\wedge\eta\rangle(z_0)\ge\lambda_{p,q}(z_0)|dz^I\wedge\eta|^2(z_0)\]
By continuity, possibly shrinking $U$, we can assume
\[\langle[\sqrt{-1}R^L,\Lambda](dz^I\wedge\eta),dz^I\wedge\eta\rangle(z)\ge(\lambda_{p,q}(z_0)-\epsilon)|dz^I\wedge\eta|^2(z)\]
holds on $U$. Applying theorem (\ref{thm:Main theorem}) to the subbundle $F=\langle dz^I\rangle$ we get the result.
\end{proof}
\end{proposition}

Next we introduce an easy lemma to glue the estimates on each cover.
\begin{proposition}\label{prop:Gluing results}
Assume there is an open covering $\mathfrak{U}=\{U_{\alpha}\}_{\alpha\in\mathscr{A}}$ of $M$ and on each $U_{\alpha}$, there exists a nonnegative function $f_{\alpha}$ on $U_{\alpha}$ such that for all $s\in\Omega_c^{0,q}(U_{\alpha},E)$,
\[\frac{1}{2}||Ds||^2\ge\int_{U_{\alpha}}f_{\alpha}|s|^2\]
Take a partition of unity $\{\chi_{\alpha}\}$ subordinate to $\mathfrak{U}$ and let $g_{\mathfrak{U}}(z)=\sum_{\alpha}|\bar{\partial}\chi_{\alpha}(z)|^2$.

Then, we have
\[\frac{1}{2}||Ds||^2\ge\int_{M}(\frac{1}{2}(\sum_{\alpha:\chi_{\alpha}(z)\neq0}f_{\alpha}^{-1})^{-1}-g_{\mathfrak{U}})|s|^2\,d\mu_g(z)\]
for all $s\in\Omega^{0,q}_c(M,E)$.
\begin{proof}
Apply the assumed estimate on $U_{\alpha}$ to $\chi_{\alpha}s$, we get
\[\int_{U_{\alpha}}|\bar{\partial}(\chi_{\alpha}s)|^2+\int_{U_{\alpha}}|\bar{\partial}^*(\chi_{\alpha}s)|^2\ge\int_{U_{\alpha}}f_{\alpha}|\chi_{\alpha}|^2|s|^2\]
Note that
\[\int_{U_{\alpha}}|\bar{\partial}(\chi_{\alpha}s)|^2=\int_{U_{\alpha}}|\bar{\partial}\chi_{\alpha}\wedge s+\chi_{\alpha}\wedge\bar{\partial}s|^2\le 2\int_{U_{\alpha}}|\mathscr{E}_{\bar{\partial}\chi_{\alpha}}s|^2+2\int_{U_{\alpha}}|\chi_{\alpha}|^2|\bar{\partial}s|^2\]
\[\int_{U_{\alpha}}|\bar{\partial}^*(\chi_{\alpha}s)|^2=\int_{U_{\alpha}}|-\partial_i\chi_{\alpha}\cdot \mathscr{I}_{dz^i}s+\chi_{\alpha}\wedge\bar{\partial}^*s|^2\le 2\int_{U_{\alpha}}|\mathscr{I}_{\partial\chi_{\alpha}}s|^2+2\int_{U_{\alpha}}|\chi_{\alpha}|^2|\bar{\partial}^*s|^2\]
And that
\[\int_{U_{\alpha}}|\mathscr{E}_{\bar{\partial}\chi_{\alpha}}s|^2+|\mathscr{I}_{\partial\chi_{\alpha}}s|^2=\int_{U_{\alpha}}\langle[\mathscr{E}_{\bar{\partial}\chi_{\alpha}},\mathscr{I}_{\partial\chi_{\alpha}}]s,s\rangle=\int_{U_{\alpha}}|\bar{\partial}\chi_{\alpha}|^2|s|^2\]
Summing and notice that $\sum_{\alpha}|\chi_{\alpha}|^2\le\sum_{\alpha}|\chi_{\alpha}|=1$, we get
\[2(\int_{M}|\bar{\partial}s|^2+\int_{M}|\bar{\partial}^*s|^2)\ge\int_{M}\sum_{\alpha}f_{\alpha}|\chi_{\alpha}|^2|s|^2-2\int_M g_{\mathfrak{U}}|s|^2\]
Finally use Cauchy-Schwarz
\[1=|\sum\chi_{\alpha}|^2\le(\sum f_{\alpha}|\chi_{\alpha}|^2)(\sum_{\alpha:\chi_{\alpha}(z)\neq0}f_{\alpha}^{-1})\]

\end{proof}

\end{proposition}

Notice that if $g_{\mathfrak{U}}(z)$ has support in $\bigcup_{\alpha\neq\beta}U_{\alpha}\cap U_{\beta}$. So if each $U_{\alpha}\cap U_{\beta}$ is precompact, $g_{\mathfrak{U}}(z)$ is bounded, and then the integrand will be bounded below by a positive constant on $\bigcup_{\alpha\neq\beta}U_{\alpha}\cap U_{\beta}$ if each $f_{\alpha}$ is bounded below by a sufficiently large constant there.

Apply gluing (\ref{prop:Gluing results}) to the previous proposition, we can state the compact case.

\begin{corollary}
If $(M,\Theta)$ is compact Hermitian manifold and $\lambda_{d-q+1}(z)>0$ everywhere or $\lambda_{d-q}(z)<0$ everywhere, then $||D_ks||^2\ge Ck||s||^2$ for all $s\in\Omega^{0,q}(M,E\otimes L^k)$ and $k>>0$. 
\end{corollary}

Gluing two type of estimates, we can state a more general corollary for the estimate.

\begin{corollary}\label{most general setting}
Let $(M,\Theta)$ be a Hermitian manifold, $E$ and $L$ are holomorphic vector bundle and line bundle with Hermitian metrics. Assume there exists a finite open covering $\{U_{\alpha}\}$ so that $U_{\alpha}\cap U_{\beta}$ is precompact for $\alpha\neq\beta$, and for every $\alpha$ there exists a $p_{\alpha}\in[d]$ with $\lambda_{p_{\alpha},q}\ge0$ on $U_{\alpha}$, $\lambda_{p_{\alpha},q}>0$ on $\overline{U_{\alpha}\cap U_{\beta}}$ for all $\beta\neq \alpha$. Moreover, assume that either one of the two following conditions holds:
\begin{itemize}
    \item $U_{\alpha}$ is precompact and $\lambda_{p_{\alpha},q}>0$ on $\bar{U}_{\alpha}$.
    \item $\lambda_{p_{\alpha}}>\lambda_{p_{\alpha}+1}$ and $\prescript{}{k}{f}_{p_\alpha,q}\ge0$ for $k>>0$ on $U_{\alpha}$.
\end{itemize}
Then for $k>>0$ there exists functions $\prescript{}{k}{f}\ge0$ such that $||D_ks||^2\ge \int_M\prescript{}{k}{f}|s|^2$ for all $s\in\Omega^{0,q}_c(M,E\otimes L^k)$. For every open set $U\subset\subset M$ with $\lambda_{p_{\alpha},q}>0$ on $\overline{U\cap U_{\alpha}}$, there exists a $C_U>0$ such that $\prescript{}{k}{f}\ge C_{U}k$.
\end{corollary}

\begin{corollary}
On $\C^d$ with metirc on trivial line bundle $L=\C$ induced by $\phi$. If
\begin{itemize}
    \item $\lambda_{d-q+1}(z)>0$ everywhere or $\lambda_{d-q}(z)<0$ everywhere
    \item $\phi=\sum_{i=1}^d\lambda_i^0|z_i|^2$ outside a compact set, $\lambda_i^0$ fixed.
\end{itemize}
Then $||D_ks||^2\ge Ck||s||^2$ for all $s\in\Omega^{0,q}_c(\C^d,L^k)$ and $k>>0$. 
\end{corollary}

Let us recall the well-known Hörmander's $L^2$-estimate and existence theorem for $\bar{\partial}$-operator \hyperlink{ref:Hör}{[Hör]}:
\begin{thm}[Hörmander-Andreotti-Vesentini]\label{L^2 estimate for bar partial}
Let $(M,\Theta)$ be a complete Hermitian manifold.
Assume that there a function $f\ge0$ such that for all $s\in\Omega^{0,q}_c(M,E)$,
\[\frac{1}{2}||Ds||^2\ge\int_Mf|s|^2\]
Then, for $b\in L^2_{(0,q)}(M,E)$ with $\bar{\partial}b=0$ in current sense and 
\[\int_M f^{-1}|b|^2\le C<\infty\]
there exists a $a\in L^2_{(0,q-1)}(M,E)$ with $\bar{\partial}a=b$ in current sense and $||a||^2\le C$.

\end{thm}

Using this, we see that in the setting of (\ref{most general setting}), if $M$ is complete, we can solve the $\bar{\partial}$-equation for $k>>0$ with the weight functions $\prescript{}{k}{f}$.

We refer to \hyperlink{ref:MM1}{[MM1]} \hyperlink{ref:Bern}{[Bern]} for the usual proof. Also, we will prove the $L^2$-estimates and existence theorems for the Dirac-Dolbeault operator $D$, which will be a generalization of the usual theorem. See (\ref{generalization of bar partial L^2 estimate}).

\subsection{Solving Dirac-Dolbeault equation}

The method of this subsection works for all choices of subspaces $\Dom\bar{\partial},\Dom\bar{\partial}^*$ satisfying assumption (\ref{domain assumption}):
\begin{assumption}\label{domain assumption}\mbox{}
\begin{itemize}
    \item $\Dom\bar{\partial},\Dom\bar{\partial}^*$ both contain $\Omega^{0,\bullet}_c(M,E)$.
    \item They can be described as Hilbert space adjoint of each other. That is,
    \[\Dom\bar{\partial}^{(q)}=\{\eta_q\in L^2_{(0,q)}(M,E):\langle\langle \bar{\partial}^*\eta_{q+1},\eta_q\rangle\rangle\le C||\eta_{q+1}||\qtext{for all}\eta_{q+1}\in\Dom\bar{\partial}^{*(q+1)}\}\]
    \[\Dom\bar{\partial}^{*(q)}=\{\eta_q\in L^2_{(0,q)}(M,E):\langle\langle \bar{\partial}\eta_{q-1},\eta_q\rangle\rangle\le C||\eta_{q-1}||\qtext{for all}\eta_{q-1}\in\Dom\bar{\partial}^{(q-1)}\}\]
    \item $\bar{\partial},\bar{\partial}^*$ sends domain back to domain:
    \[\im\bar{\partial}\subset\ker\bar{\partial}\]
    \[\im\bar{\partial}^*\subset\ker\bar{\partial}^*\]
\end{itemize}
\end{assumption}

We then set the domain of $D=\sqrt{2}(\bar{\partial}+\bar{\partial}^*)$ to be
\[\Dom D=\Dom\bar{\partial}\cap\Dom\bar{\partial}^*\]
Finally define
\[\Dom\Box=\{s\in\Dom D:Ds\in\Dom D\}=\{s\in\Dom\bar{\partial}\cap\Dom\bar{\partial}^*:\bar{\partial}s\in\Dom\bar{\partial}^*,\bar{\partial}^*s\in\Dom\bar{\partial}\}\]
We get a self-adjoint extension of the Kodaira Laplacian $\Box$. Notice that if $s\in\ker\Box$, then
\[||Ds||^2=2\langle\langle\Box s,s\rangle\rangle=0\]
Thus $\ker\Box=\ker D$, and the Bergman projection $B_{\Box}:L^2_{(0,\bullet)}(M,E)\to\ker\Box$ is exactly the orthogonal projection $B_{D}:L^2_{(0,\bullet)}(M,E)\to\ker D$.

$\ker\bar{\partial}$ and $\ker\bar{\partial}^*$ are both closed and for any $\eta\in L^2(M,E)$ the orthogonal projections $\eta^{\ker\bar{\partial}},\eta^{\ker\bar{\partial}^*}$ satisfy
\[\eta-\eta^{\ker\bar{\partial}}\in \ker\bar{\partial}^*\]
\[\eta-\eta^{\ker\bar{\partial}^*}\in\ker\bar{\partial}\]
Thus, these two projections send $\Dom\bar{\partial}$ and $\Dom\bar{\partial}^*$ back into themselves. For example if $\eta\in\Dom D$, then $\eta^{\ker\bar{\partial}},\eta^{\ker\bar{\partial}^*}\in\Dom D$.

The following is a generalization of the technique in \hyperlink{ref:Hör}{[Hör]}.

\begin{thm}\label{thm:L2 estimate for Hodge-Dolbeault}
Assume there exists nonnegative measurable functions $f_{q+1},f_{q-1}\ge0$ on $M$ such that for all $\eta_{q+1}\in\ker\bar{\partial}^{(q+1)}\cap\Dom \bar{\partial}^{*(q+1)}$,
\[||\bar{\partial}^*\eta_{q+1}||^2\ge\int_Mf_{q+1}|\eta_{q+1}|^2\]
and for $\eta_{q-1}\in\ker\bar{\partial}^{*(q-1)}\cap\Dom \bar{\partial}^{(q-1)}$,
\[||\bar{\partial}\eta_{q-1}||^2\ge\int_{M}f_{q-1}|\eta_{q-1}|^2\]
Then we have the following $L^2$-estimate for the Dirac-Dolbeault operator $D^{(q)}$ on $(0,q)$-forms:

For any $b_{q+1}\in\ker\bar{\partial}^{(q+1)}$ and $b_{q-1}\in\ker\bar{\partial}^{*(q-1)}$ with
\[\int_M f_{q+1}^{-1}|b_{q+1}|^2+f_{q-1}^{-1}|b_{q-1}|^2\le C<\infty\]
there exists a $a_q\in \Dom D^{(q)}$ so that $D^{(q)}a_q=b_{q+1}+b_{q-1}$ with the estimate
\[||a_q||^2\le C\]
\begin{proof}
Define a subspace of $L^2_{(0,q)}(M,E)$
\[V=\{\bar{\partial}^*\eta_{q+1}+\bar{\partial}\eta_{q-1}:\eta_{q+1}\in\Dom\bar{\partial}^{(q+1)*},\ \eta_{q-1}\in\Dom \bar{\partial}^{(q-1)}\}\]
We define a functional $T$ on $V$ by
\[T(\bar{\partial}^*\eta_{q+1}+\bar{\partial}\eta_{q-1})=\langle\langle \eta_{q+1},b_{q+1}\rangle\rangle+\langle\langle \eta_{q-1},b_{q-1}\rangle\rangle\]
We need to check the well-definedness of this functional. If $\bar{\partial}^*\eta_{q+1}+\bar{\partial}\eta_{q-1}=0$, then $\bar{\partial}\eta_{q-1}\in\ker\bar{\partial}\cap\im\bar{\partial}^*=0$. Thus also $\bar{\partial}^*\eta_{q+1}=0$. Next, we show that $\langle\langle \eta_{q+1},b_{q+1}\rangle\rangle=0$. If $f_{q+1}=0$ a.e., then by assumption $b_{q+1}=0$ a.e. and $\langle\langle \eta_{q+1},b_{q+1}\rangle\rangle=0$. Otherwise $f_{q+1}>0$ on some set of positive measure.

We always have
\[\eta_{q+1}-\eta_{q+1}^{\ker\bar{\partial}}\in\ker\bar{\partial}^*\]
Thus also $\bar{\partial}^*\eta_{q+1}^{\ker\bar\partial}=0$. Now $\Box\eta_{q+1}=DD\eta_{q+1}=0$, by elliptic regularity, it is smooth. We apply the assumed estimate to $\eta_{q+1}^{\ker\bar{\partial}}$, we get $\eta_{q+1}=0$ on a set of positive measure. By Lebesgue density theorem, there exists a point $x$ with density 1 with respect to the zero set of $\eta_{q+1}$, then every degree of Taylor expansion at $x$ vanishes. By unique continuation of second order elliptic partial differential equation \hyperlink{ref:Aron}{[Aron]}, $\eta_{q+1}^{\ker\bar{\partial}}\equiv0$. Finally from $b_{q+1}\in\ker\bar{\partial}$,
\[\langle\langle \eta_{q+1},b_{q+1}\rangle\rangle=\langle\langle\eta_{q+1}^{\ker\bar{\partial}},b_{q+1}\rangle\rangle=0\]
Similar argument shows that either $b_{q-1}=0$ a.e. or $\eta_{q-1}^{\ker\bar{\partial}^*}=0$. We conclude that $T$ is well-defined.

Moreover, the operator norm of $T$ is bounded by $C^{1/2}$ on $V$:
\begin{equation*}
\begin{split}
&|\langle\langle \eta_{q+1},b_{q+1}\rangle\rangle+\langle\langle \eta_{q-1},b_{q-1}\rangle\rangle|\\
&=|\langle\langle \eta_{q+1}^{\ker\bar{\partial}},b_{q+1}\rangle\rangle+\langle\langle \eta_{q-1}^{\ker\bar{\partial}^*},b_{q-1}\rangle\rangle|\\
&\le(\int_Mf_{q+1}^{-1}|b_{q+1}|^2)^{1/2}(\int_{M}f_{q+1}|\eta_{q+1}^{\ker\bar{\partial}}|^2)^{1/2}+(\int_Mf_{q-1}^{-1}|b_{q-1}|^2)^{1/2}(\int_{M}f_{q-1}|\eta_{q-1}^{\ker\bar{\partial}^*}|^2)^{1/2}\\
&\le C^{1/2}(\int_{M}|\bar{\partial}^*\eta_{q+1}^{\ker\bar{\partial}}|^2+\int_{M}|\bar{\partial}\eta_{q-1}^{\ker\bar{\partial}^*}|^2)^{1/2}\\
&= C^{1/2}(\int_{M}|\bar{\partial}^*\eta_{q+1}|^2+\int_{M}|\bar{\partial}\eta_{q-1}|^2)^{1/2}\\
&=C^{1/2}||\bar{\partial}^*\eta_{q+1}+\bar{\partial}\eta_{q-1}||
\end{split}
\end{equation*}
By Hahn-Banach theorem, $T$ extends to a bounded linear functional on whole $L^2_{(0,q)}(M,E)$ with norm bounded by $C^{1/2}$. By Hilbert space theory, $T$ is represented by the inner product with an element $a_q\in L^2_{(0,q)}(M,E)$ with $||a_q||^2\le C$. Notice that then
\[\langle\langle \eta_{q+1},\bar{\partial}a_q\rangle\rangle+\langle\langle \eta_{q-1},\bar{\partial}^*a_q\rangle\rangle=\langle\langle\bar{\partial}^*\eta_{q+1}+\bar{\partial}\eta_{q-1},a_{q}\rangle\rangle=\langle\langle \eta_{q+1},b_{q+1}\rangle\rangle+\langle\langle \eta_{q-1},b_{q-1}\rangle\rangle\]
for all $\eta_{q+1}\in\Dom\bar{\partial}^*,\ \eta_{q-1}\in\Dom\bar{\partial}$, which says exactly that $D^{(q)}a_{q}=b_{q+1}+b_{q-1}$.
\end{proof}
\end{thm}

\begin{remark}\label{generalization of bar partial L^2 estimate}
We can always choose $f_{q-1}\equiv0$, since then the condition for $\eta_{q-1}$ is trivially satisfied. Then the form $b_{q-1}$ is forced to be $0$, as seen in the proof of (\ref{thm:L2 estimate for Hodge-Dolbeault}). In this special case, the theorem says that if we have the estimation for $(0,q+1)$-form, then for any $b_{q+1}\in\ker\bar{\partial}^{(q+1)}$ we can solve the equation
\[\bar{\partial}a_q=b_{q+1}\qquad\bar{\partial}^*a_q=0\]
with the estimate
\[||a_q||^2\le\int_Mf_{q+1}^{-1}|b_{q+1}|^2\]
This special case is the original $L^2$-estimate and existence theorem for $\bar{\partial}$-operator. The part $\bar{\partial}^*a_q=0$ can be obtained from considering $a_q-a_q^{\ker\bar{\partial}}$.
\end{remark}

By this theorem, we can deduce the local spectral gap for Bergman kernel.
\begin{corollary}
Assume the estimate in (\ref{thm:L2 estimate for Hodge-Dolbeault}), and let $U\subset M$ be an open set so that $f_{q+1},f_{q-1}\ge C$ on $U$. Denote $B=B_{\Box}=B_{D}$ to be the Bergman projection. Then, for $s\in\Omega^{0,q}_c(U,E)$,
\[||s-Bs||^2\le \frac{1}{2C}||Ds||^2\le \frac{1}{C^2}||\Box s||^2\]
\begin{proof}
Now $s\in\Omega^{0,q}_c(U,E)\subset\Dom D$, and we have
\[\bar{\partial}\bar{\partial}s=0\qquad\bar{\partial}^*\bar{\partial}^*s=0\]
Thus we can solve the equation $D^{(q)}a=D^{(q)}s$ with
\[||a||^2\le\int_Mf_{q+1}^{-1}|\bar{\partial}s|^2+f_{q-1}^{-1}|\bar{\partial}^*s|^2\le \frac{1}{2C}||Ds||^2\]
By the definition of $B_{D}$, $s-B_{D}s$ is a solution to $D^{(q)}a=D^{(q)}s$ with the smallest norm, we deduce
\[||s-B_{D}s||^2\le\frac{1}{2C}||Ds||^2\]
Finally, since $\bar{\partial}s\in\ker\bar{\partial}^{(q+1)}\cap\Dom\bar{\partial}^{*(q+1)}$, we can apply the estimate to $\bar{\partial}$:
\[||\bar{\partial}^*\bar{\partial}s||^2\ge\int_Mf_{q+1}|\bar{\partial}s|^2\ge C||\bar{\partial}s||^2\]
Similar for $\bar{\partial}^*s$, we get
\[\frac{1}{2C}||Ds||^2\le\frac{1}{C^2}(||\bar{\partial}^*\bar{\partial}s||^2+||\bar{\partial}\bar{\partial}^*s||^2)=\frac{1}{C^2}||\Box s||^2\]
\end{proof}
\end{corollary}

In the following section, we will use the theorem in previous section to define appropriate domains for different geometric conditions and deduce the estimation for both $(0,q-1)$ and $(0,q+1)$-forms.

\section{Aymptotics of Bergman kernel}\label{section:Bergman kernel}

In this section, we show that the expansion in local spectral gap case is equivalent to the model case (\ref{model case}) near nondegenerate points, and the model case admits a global spectral gap. The main tool is the finite propagation speed for wave equation used in \hyperlink{ref:MM1}{[MM1]}.

Let $M(q)$ be the points where $R^L$ is nondegenerate with signature $(q,d-q)$ and let $U\subset M(q)$.

We only need to assume the following assumptions (\ref{local spectral gap assumption}):
\begin{assumption}\label{local spectral gap assumption}\mbox{}
\begin{itemize}
    \item The self-adjoint extensions of $D_k$ and $\Box_k$ satisfy $\ker D_k=\ker\Box_k$.
    \item We have local spectral gap of growth order $k^{-\alpha}$ on $U$ for some $\alpha\in\R$:
    \[||s-B_{\Box_k}s||\le C_Uk^{\alpha}||\Box_ks||\]
    for all $s\in\Omega^{0,q}_c(U,E\otimes L^k)$.
\end{itemize}
\end{assumption}
Under this assumption we simply denote $B_k=B_{\Box_k}=B_{D_k}$ to be the Bergman projection.

\subsection{Complete case}

In the case that $(M,\Theta)$ is complete, several self-adjoint extensions of Kodaira Laplacian are in fact equivalent. We can first define the usual Gaffney extension as follows: Let $\bar{\partial}$ be the maximal extension, which means that $\Dom\bar{\partial}$ consists of elements $s\in L^2_{(0,\bullet)}(M,E)$ that has an $L^2$-distributional derivative: $\bar{\partial}s\in L^2_{(0,\bullet)}(M,E)$. Then we define $\bar{\partial}^*=\bar{\partial}^*_H$ to be its Hilbert space adjoint,
\[\Dom\bar{\partial}^*=\{s\in L^2_{(0,\bullet)}(M,E):\langle \bar{\partial}\eta,s\rangle\le C|\eta|\qtext{for all}\eta\in\Dom\bar{\partial}\}\]
This satisfies (\ref{domain assumption}), and we set the domain of $D,\Box$ as the discussion after the assumption.

Due to completeness, one can prove that $\bar{\partial}^*$ is actually the maximal extension \hyperlink{ref:MM1}{[MM1, Cor.3.3.3]}. If we define a positive quadratic form $Q$ on test forms:
\[Q(s_1,s_2)=\langle\langle\bar{\partial}s_1,\bar{\partial}s_2\rangle\rangle+\langle\langle\bar{\partial}^*s_1,\bar{\partial}^*s_2\rangle\rangle\]
Then we can prove that test forms are dense in $\Dom D$ with the norm
\[||s||_Q^2=||s||^2_{L^2}+Q(s,s)\]
Thus $\Box$ coincides with the Friedrichs extension. As a result, our estimate holds for all $s\in\Dom D^{(q)}$:
\[\frac{1}{2}||D_ks||^2\ge\int_M\prescript{}{k}{f}_{p,q}|s|^2\]

Combining all the discussions before, we can formulate a condition on complete manifolds.
\begin{corollary}[Local Spectral gap]
Let $(M,\Theta)$ be a complete Hermitian manifold. $L$ and $E$ are holomorphic line bundle and vector bundle over $M$ with Hermitian metrics. Let $\lambda_1(z)\ge\cdots\ge\lambda_d(z)$ be the eigenvalues of $R^L(z)$, and assume the following conditions:
\begin{itemize}
    \item $\lambda_{d-q}\ge0\ge\lambda_{d-q+1}$ and $\lambda_{d-q}>\lambda_{d-q+1}$ everywhere.
    \item The functions $\prescript{}{k}{f}_{d-q,q+1},\prescript{}{k}{f}_{d-q,q-1}\ge0$ for $k>>0$. 
    
    In other words, the function $\lambda^{E}_{d-q}+C_d|\partial\Theta|^2+|\bar{\partial}\Psi^{F_{d-q}}|^2$ is bounded by $k\lambda_{d-q}$ and $-k\lambda_{d-q+1}$ on $M$ for $k>>0$.
\end{itemize}
Then, for any open set $U\subset\subset M$ with $\lambda_{d-q}(z)>0>\lambda_{d-q+1}(z)$ everywhere on $\overline{U}$, there exists a constant $C_U$ so that for all $s\in\Omega^{0,q}_c(U,E\otimes L^k)$ and $k>>0$,
\[||s-B_ks||^2_k\le\frac{C_U}{2k}||D_ks||_k^2\le\frac{C_U^2}{k^2}||\Box_ks||_k^2\]
That is, there is a local spectral gap of growth order $k$ on $U$, (\ref{local spectral gap assumption}) is satisfied.

\end{corollary}

\begin{example}
Let $\phi_1,\phi_2$ be real functions defined on $\C^{d-q}$ and $\C^{q}$ respectively with
\[\sqrt{-1}\partial\bar{\partial}\phi_1\ge0\qquad\sqrt{-1}\partial\bar{\partial}\phi_2\le0\]
Then for the trivial line bundle $L=\C$ over $\C^d=\C^{d-q}\times \C^q$ with metric defined by $e^{-\phi_1-\phi_2}$, then
\[\partial\bar{\partial}\phi=\sum_{i,j\le d-q}(\phi_1)_{i\bar{j}}\,dz^i\wedge d\bar{z}^j+\sum_{i,j\ge d-q+1}(\phi_2)_{i\bar{j}}\,dz^i\wedge d\bar{z}^j\]

We pick $F=\langle dz^1\wedge\cdots\wedge dz^{d-q}\rangle$ (we do not need $\lambda_{d-q}>\lambda_{d-q+1}$ here) and notice that
\[\Psi^F=dz^1\wedge\cdots\wedge dz^{d-q}\otimes \partial_1\wedge\cdots\wedge\partial_{d-q}\]
Thus $\bar{\partial}\Psi^F=0$. As a result, for $L^k$-valued forms,
\[\prescript{}{k}{f}_{d-q,q+1}=k\lambda_{d-q}\ge0\qquad\prescript{}{k}{f}_{d-q,q-1}=-k\lambda_{d-q+1}\ge0\]
We obtain local spectral gap of growth order $k$ (\ref{local spectral gap assumption}) on $U\subset\subset M(q)$.

For example, if we take Schwartz functions $\alpha^i:\C\to\R$ with $\alpha^i\ge0$ for $i\le d-q$ and $\alpha^i\le0$ for $i\ge d-q+1$, and set
\[\partial\bar{\partial}\phi=\sum_{i=1}^d\alpha^i(z^i)\,dz^i\wedge d\bar{z}^i\]
Then we can pick $\phi$ to be the sum of Newtonian potentials associated to $\alpha^i$. In this example $|\alpha^i|$ have no uniform lower bound and we even can pick compactly supported ones.
\end{example}

\begin{example}\label{example:complete product complete}
More generally, for $(M_1,\Theta_1,L_1),(M_2,\Theta_2,L_2)$ both complete of dimension $d-q,q$, assume
\[\sqrt{-1}R^{L_1}\ge0\qquad\sqrt{-1}R^{L_2}\le0\]
Consider the product $M_1\times M_2$, and let $p_1,p_2$ be two projections and $L=p_1^*L_1\otimes p_2^*L_2$. Locally take $\{\partial_i\}_{i=1}^{d-q}$ a holomorphic coordinate for $M_1$, we still have
\[\Psi^F=dz^1\wedge\cdots\wedge dz^{d-q}\otimes \partial_1\wedge\cdots\wedge\partial_{d-q}\]
and $\bar{\partial}\Psi^F=0$. So for $E\to M_1\times M_2$,
\[\prescript{}{k}{f}_{d-q,q+1}=\frac{2}{3}k\lambda_{d-q}-\frac{2}{3}\lambda^E_{d-q}-C_d|\partial\Theta|^2\qquad\prescript{}{k}{f}_{d-q,q-1}=-\frac{2}{3}k\lambda_{d-q+1}-\frac{2}{3}\lambda^E_{d-q}-C_d|\partial\Theta|^2\]
If these two functions are both $\ge0$, we have local spectral gap of growth order $k$ (\ref{local spectral gap assumption}) on $U\subset\subset M(q)$.
\end{example}

We mention a case that attains a global spectral gap, and in the last subsection we will show that the asymptotic behavior of Bergman kernel at nondegenerate points satisfying (\ref{local spectral gap assumption}) is equivalent to the following case.

\begin{example}[model case]\label{model case}
Now on $\C^d$ with $L=\C,E=\C^r$ trivial bundles, and the metrics are given by $(g,e^{-\phi},h)$. If $\partial\bar{\partial}\phi$ is nondegenerate with fixed signature $(q,d-q)$ over all $\C^d$ and outside a compact set we have
\begin{itemize}
    \item $g=g_0$, where $(g_0)_{i\bar{j}}=\delta_{ij}$.
    \item $h=h_0$, where $(h_0)_{a\bar{b}}=\delta_{ab}$ for global trivializing sections $\{e_a\}_{a=1}^r$.
    \item $\phi=\phi_0$, where $\phi_0=\sum_{i=1}^d\lambda_{i}^0|z^i|^2$ with
    \[\lambda_1^0\ge\cdots\ge\lambda_{d-q}^0>0>\lambda_{d-q+1}^0\ge\cdots\ge\lambda_{d}^0\]
\end{itemize}
Then we have $\Spec D_k\subset\{0\}\cup[C\sqrt{k},\infty)$ for some $C>0$ and $k>>0$.
\end{example}

Next, we take a closer look into the case that $\dim_{\C}M=2$ and the signature of $R^L$ is $(1,1)$, which is the smallest nontrivial mixed curvature case.

\begin{proposition}\label{dimension 2}
Let $(M,\Theta)$ be a $2$-dimensional complete Hermitian manifold. $E,L$ are holomorphic vector bundle and line bundles over $M$ with metrics and we assume that $R^L$ has fixed signature $(1,1)$. If there exists four positive constant $C_1,C_2,C_3,C_4$ such that
\begin{itemize}
    \item $\det R\le -C_1$
    \item $[\sqrt{-1}R^{E\otimes T^{1,0}},\Lambda]\ge -\frac{C_2}{|\tr R|+1}\id^{E\otimes T^{1,0}}$
    \item $|\partial\Theta|^2\le\frac{C_3}{|\tr R|+1}$
    \item $|\bar{\partial}R|^2\le C_4(|\tr R|+1)$
\end{itemize}
Then we have local spectral gap of growth order $k$ (\ref{local spectral gap assumption}) for $(0,1)$-forms on every $U\subset\subset M(q)$.
\begin{proof}
For two dimensional case, the eigenvalues are easy to write, for example
\[\lambda_{1,2}=\lambda_1=\frac{1}{2}(\tr R+\sqrt{(\tr R)^2-4\det R})\ge\frac{1}{2}(\tr R+\sqrt{(\tr R)^2+4C_1})\]
\[\lambda_{1,0}=-\lambda_2=\frac{1}{2}(-\tr R+\sqrt{(\tr R)^2-4\det R})\ge\frac{1}{2}(-\tr R+\sqrt{(\tr R)^2+4C_1})\]
\[(\lambda_2-\lambda_1)^2=(\tr R)^2-4\det R\ge(\tr R)^2+4C_1\]
and we want for $k>>0$,
\[\frac{2}{3}k\lambda_{1,2}-\lambda^E_p-C_d|\partial\Theta|^2-\frac{|\bar{\partial}R|^2}{(\lambda_2-\lambda_1)^2}\ge0\]
By assumption, as $\tr R\to-\infty$, 
\[\lambda_{1,2}\ge\frac{C}{|\tr R|+1},\qquad \lambda_p^E+C_d|\partial\Theta|^2+\frac{|\partial R|^2}{(\lambda_2-\lambda_1)^2}\le \frac{C^\prime}{|\tr R|+1}\]
Thus the inequality holds for $k>>0$ and $\tr R<<0$. The part $\tr R\ge C_0$ is easy.
Similarly for $(0,2)$-form we have the corresponding argument. Apply previous theorem and we are done.

\end{proof}
\end{proposition}

\begin{example}\label{key example}
We construct an explicit example on $\C^2$. Let $f:\C^2\to\C^2$ be holomorphic. Write $f=(f^1,f^2)$. Consider the closed form 
\[\omega=\alpha^1(z^1)\,dz^1\wedge d\bar{z}^1+\alpha^2(z^2)\,dz^2\wedge d\bar{z}^2\]
where $\alpha^1\ge0,\ \alpha^2\le0$ are real-valued function on $\C$. Consider the pullback $f^*\omega$, then there exists real function $\phi$ on $\C^2$ such that
\[\partial\bar{\partial}\phi=f^*\omega\]
Let $e^{-\phi}$ define the metric of trivial line bundle $L=\C$.

Compute
\begin{equation*}
\begin{split}
\matrixtwo{R_{1\bar{1}}}{R_{2\bar{1}}}{R_{1\bar{2}}}{R_{2\bar{2}}}
&=\overline{J_f}^t\matrixtwo{\alpha^1}{0}{0}{\alpha^2}J_f\\
&=\matrixtwo{\bar{f^1_1}}{\bar{f^2_1}}{\bar{f^1_2}}{\bar{f^2_2}}\matrixtwo{\alpha^1}{0}{0}{\alpha^2}\matrixtwo{f^1_1}{f^1_2}{f^2_1}{f^2_2}\\
&=\matrixtwo{\alpha^1|f^1_1|^2+\alpha^2|f^2_1|^2}{\alpha^1\bar{f^1_1}f^1_2+\alpha^2\bar{f^2_1}f^2_2}{\alpha^1f^1_1\bar{f^1_2}+\alpha^2f^2_1\bar{f^2_2}}{\alpha^1|f^1_2|^2+\alpha^2|f^2_2|^2}
\end{split}
\end{equation*}
\begin{equation*}
\begin{split}
R=\matrixtwo{R^1_1}{R^1_2}{R^2_1}{R^2_2}=\matrixtwo{g^{1\bar{1}}}{g^{1\bar{2}}}{g^{2\bar{1}}}{g^{2\bar{2}}}\matrixtwo{\bar{f^1_1}}{\bar{f^2_1}}{\bar{f^1_2}}{\bar{f^2_2}}\matrixtwo{\alpha^1}{0}{0}{\alpha^2}\matrixtwo{f^1_1}{f^1_2}{f^2_1}{f^2_2}
\end{split}
\end{equation*}

If $g_{i\bar{j}}=u^{-1}\delta_{ij}$ for some positive function $u$,
\[\lambda_1\cdot\lambda_2=u^{2}\alpha^1\alpha^2|f^1_1f^2_2-f^1_2f^2_1|^2\]
\[\lambda_1+\lambda_2=u[\alpha^1(|f^1_1|^2+|f^1_2|^2)+\alpha^2(|f^2_1|^2+|f^2_2|^2)]\]
\[\lambda_1-\lambda_2=u[(\alpha^1(|f^1_1|^2+|f^1_2|^2)+\alpha^2(|f^2_1|^2+|f^2_2|^2))^2-4\alpha^1\alpha^2|f^1_1f^2_2-f^1_2f^2_1|^2]^{1/2}\]

We further take $\alpha^1=1,\alpha^2=-1$, and $f^1=z^1+p$, $f^2=z^2+p$, then
\[\lambda_1\cdot\lambda_2=-u^2|1+p_1+p_2|^2\]
We require $1+p_1+p_2\neq0$ everywhere on $\C^2$. (For example $p=\sum a_l(z^1-z^2)^l$, then $p_1+p_2=0$.)
\[\lambda_1+\lambda_2=2u\Re (p_1-p_2)\]
\[\lambda_1-\lambda_2=2u[(\Re(p_1-p_2))^2+|1+p_1+p_2|^2]^{1/2}\]
\[\lambda_{1,2}=\lambda_1=u\{\Re(p_1-p_2)+[(\Re(p_1-p_2))^2+|1+p_1+p_2|^2]^{1/2}\}\]
\[\lambda_{1,0}=-\lambda_2=u\{-\Re(p_1-p_2)+[(\Re(p_1-p_2))^2+|1+p_1+p_2|^2]^{1/2}\}\]
\[R=u\matrixtwo{1+2\Re p_1}{p_2-\bar{p_1}}{\bar{p_2}-p_1}{-1-2\Re p_2}\]
Thus,
\begin{equation*}
\begin{split}
|\bar{\partial}R|^2
&=u\{|u_{\bar{1}}(1+2\Re p_1)+u\bar{p_{11}}|^2+|u_{\bar{1}}(p_2-\bar{p_1})-u\bar{p_{11}}|^2+|u_{\bar{1}}(\bar{p_2}-p_1)+u\bar{p_{21}}|^2+|u_{\bar{1}}(-1-2\Re p_2)-u\bar{p_{21}}|^2\\
&+|u_{\bar{2}}(1+2\Re p_1)+u\bar{p_{12}}|^2+|u_{\bar{2}}(p_2-\bar{p_1})-u\bar{p_{12}}|^2+|u_{\bar{2}}(\bar{p_2}-p_1)+u\bar{p_{22}}|^2+|u_{\bar{2}}(-1-2\Re p_2)-u\bar{p_{22}}|^2\}
\end{split}
\end{equation*}
Now if $u\equiv 1$ then $|\bar{\partial}R|^2=2\sum_{i,j=1,2}|p_{ij}|^2$. Consider $p=(z^1-z^2)^2$, then $|\bar{\partial}R|^2=32$ and $|\tr R|=4\Re(z^1-z^2)$. The condition in (\ref{dimension 2}) is satisfied, and for any $U\subset\subset\C^n$ we have local spectral gap of growth order $k$ (\ref{local spectral gap assumption}). For the explicit metric $\phi$ we can pick 
\begin{equation*}
\begin{split}
\phi
&=|z^1|^2-|z^2|^2\\
&+(z^1)^2\bar{z}^1+z^1(\bar{z}^1)^2-(z^2)^2\bar{z}^2-z^2(\bar{z}^2)^2+z^1(\bar{z}^2)^2+\bar{z}^1(z^2)^2-(z^1)^2\bar{z}^2-(\bar{z}^1)^2z^2\\
&-2z^1\bar{z}^1z^2-2z^1\bar{z}^1\bar{z}^2+2z^1z^2\bar{z}^2+2\bar{z}^1z^2\bar{z}^2
\end{split}
\end{equation*}
In this example, $\lambda_1,-\lambda_2$ has no uniform lower bound, the error $|\bar{\partial}\Psi^{F_1}|^2$ is not zero, and the bound $(\lambda_1-\lambda_2)^{-2}|\bar{\partial}R|^2$ suffices.

\end{example}

\subsection{Neumann boundary condition}

In this subsection, we consider Hermitian manifold with smooth boundary $\partial M$. For any $z\in\partial M$, we can locally take a smooth real-valued function $\rho$ defined on a coordinate neighborhood such that the interior of $M$ is described by $\rho<0$ and $|d\rho|=1$ on $\partial M$. Then we can describe the outward unit normal 
$\mathfrak{n}$ as the metric dual $(d\rho)^*$. Given a frame $\{w_i\}$ of $T^{1,0}M$, we can write
\[\mathfrak{n}^{1,0}=(\partial\rho)^*=\rho_{\bar{j}}g^{i\bar{j}}w_i\qquad \mathfrak{n}^{0,1}=(\bar{\partial}\rho)^*=\rho_ig^{i\bar{j}}\bar{w}_j\]

Given a line bundle $L\to M$, for any $z\in\overline{M}$ we can unitarily diagonalized $R^L$, so that $T_zM=E_+\oplus E_-$, $E_+$ is the orthogonal sum of eigenspaces of $\lambda_i$, $i=1,\dots,p$ and $E_-$ is the direct sum of eigenspaces of $\lambda_j$, $j=p+1,\dots,d$. We use $+,-$ to indicate that for $R^L$ nondegenerate of signature $(d-p,p)$, $\lambda_p>0>\lambda_{p+1}$ and $E_+,E_-$ represents the direct sum of positive, negative eigenspaces. However, for most of the discussions in this subsection we do not need to assume the sign of eigenvalue.

In the subsection, we consider two cases: Either $\mathfrak{n}^{1,0}\in E_+$ along $\partial M$, or $\mathfrak{n}^{1,0}\in E_-$ along $\partial M$. 

Notice that if $p=d$, the assumption $\mathfrak{n}^{1,0}\in E_+$ is trivially satisfied, which corresponds to the $\lambda_i\ge0$ case that we have no restriction on normal directions. Similarly $p=0$ corresponds to the $\lambda_i\le0$ case with no extra normal assumptions.

Now assume $\lambda_p>\lambda_{p+1}$, then we can pick a smooth local frame $\{w_i\}$ with index set $[+]\cup[-]=[d]$ such that $w_i\in E_+$ for $i\in[+]$ and $w_j\in E_-$ for $j\in[-]$. Then if $\mathfrak{n}^{1,0}\in E_+$, we have $\rho_{\bar{j}}=0=\rho_j$ along $\partial M$ for $j\in[-]$.

If $\mathfrak{n}^{1,0}\in E_-$, we have $\rho_i=0=\rho_{\bar{i}}$ along $\partial M$ for $i\in[+]$.

We will use introduce two different boundary conditions and extensions of $D$. First for the case $\mathfrak{n}^{1,0}\in E_+$ we use the usual Gaffney extension as described in the complete case. $\bar{\partial}$ is the maximal extension and $\bar{\partial}^*$ is its Hilbert space adjoint. However, if $M$ has boundary, $\bar{\partial}^*$ will not be the maximal extension. In fact, for $s_1,s_2\in\Omega^{\bullet,\bullet}_{c}(\bar{M},E\otimes L^k)$,
\[\langle\langle\bar{\partial}s_1,s_2\rangle\rangle-\langle\langle s_1,\bar{\partial}^*s_2\rangle\rangle=\int_{\partial M}\langle s_1,\mathscr{I}_{\partial\rho}s_2\rangle\]
From this we see $\Omega^{\bullet,\bullet}_{c}(\bar{M},E\otimes L^k)\subset\Dom\bar{\partial}$ and denote
\[B^{\bullet,\bullet}_{+}:=\Dom\bar{\partial}^*\cap \Omega^{\bullet,\bullet}_{c}(\bar{M},E\otimes L^k)=\{s\in\Omega^{\bullet,\bullet}_{c}(\bar{M},E\otimes L^k):\mathscr{I}_{\partial\rho}s=0\}\]
As in \hyperlink{ref:Bern}{[Bern]}, using the well-know Friedrichs lemma carefully along the boundary, we know $B_+^{0,q}$ is dense in $\Dom\bar{\partial}\cap\Dom\bar{\partial}^*$ with respect to the norm
\[||s||^2_Q=||s||^2_{L^2}+||\bar{\partial}s||^2_{L^2}+||\bar{\partial}^*s||^2_{L^2}\]

Next if $\mathfrak{n}^{0,1}\in E_-$, we do the dual process as following. Let $\bar{\partial}^*$ be the maximal extension and $\bar{\partial}$ be its Hilbert space adjoint. This time $\Omega^{\bullet,\bullet}_{c}(\bar{M},E\otimes L^k)\subset\Dom\bar{\partial}^*$ and
\[B^{\bullet,\bullet}_{-}:=\Dom\bar{\partial}\cap \Omega^{\bullet,\bullet}_{c}(\bar{M},E\otimes L^k)=\{s\in\Omega^{\bullet,\bullet}_{c}(\bar{M},E\otimes L^k):\mathscr{E}_{\bar{\partial}\rho}s=0\}\]
Similarly $B_-^{0,q}$ is dense in $\Dom\bar{\partial}\cap\Dom\bar{\partial}^*$ with respect to $||\cdot||_Q$.

The domains in two cases are both defined by Gaffney's extension process and satisfies (\ref{domain assumption}).

We need to derive the estimate with boundary term under mixed boundary condition.

\begin{thm}\label{L^2 estimates with boundary}
Assume $\lambda_p>\lambda_{p+1}$, $F=F_p$ as in (\ref{thm:Main theorem applied to important bundle})
\begin{enumerate}[(a)]
    \item If $\mathfrak{n}^{1,0}\in E_+$ and $s\in B^{0,q}_+$,
    \[\frac{1}{2}||D_ks||^2\ge\int_M \prescript{}{k}{f}_{p,q}|s|^2+\int_{\partial M}\langle(\partial\bar{\partial}\rho)(w_i,\bar{w}_j)\mathscr{E}_{\bar{w}^j}\mathscr{I}_{w^i}s,s\rangle-\langle\mathscr{E}_{\partial\rho}\mathscr{I}_{\bar{w}^j}\tilde{\nabla}_{\bar{j}}\Psi^F,\Psi^F\rangle|s|^2\]
    \item $\mathfrak{n}^{1,0}\in E_-$ and $s\in B^{0,q}_{-}$,
    \[\frac{1}{2}||D_ks||^2\ge\int_M \prescript{}{k}{f}_{p,q}|s|^2+\int_{\partial M}\langle(\partial\bar{\partial}\rho)(w_i,\bar{w}_j)\mathscr{I}_{w^i}\mathscr{E}_{\bar{w}^j}s,s\rangle-\langle\mathscr{I}_{\bar{\partial}\rho}\mathscr{E}_{w^i}\tilde{\nabla}_i\Psi^F,\Psi^F\rangle|s|^2\]
\end{enumerate}
\begin{proof}
Apply the Bochner-Kodaira-Nakano inequality with boundary terms (\ref{BKN inequality with boundary terms}):
\begin{equation*}
\begin{split}
\frac{3}{2}||D_k(\Psi^F\wedge s)||^2
&\ge\langle\langle[\sqrt{-1}R^{E\otimes\bigwedge\nolimits_{p,0}\otimes L^k},\Lambda]\Psi^F\wedge s,\Psi^F\wedge s\rangle\rangle\\
&-||\mathcal{T}(\Psi^F\wedge s)||^2-||\mathcal{T}^*(\Psi^F\wedge s)||^2-||\bar{\mathcal{T}}(\Psi^F\wedge s)||^2-||\bar{\mathcal{T}}^*(\Psi^F\wedge s)||^2\\
&+\int_{\partial M}\langle\bar{\partial}(\Psi^F\wedge s),\bar{\partial}\rho\wedge \Psi^F\wedge s\rangle+\int_{\partial M}\langle\partial\rho\wedge(\nabla^{1,0*}+\mathcal{T}^*)(\Psi^F\wedge s),\Psi^F\wedge s\rangle\\
&-\int_{\partial M}\langle\bar{\partial}\rho\wedge(\bar{\partial}^*+\bar{\mathcal{T}}^*)(\Psi^F\wedge s),\Psi^F\wedge s\rangle-\int_{\partial M}\langle\nabla^{1,0}(\Psi^F\wedge s),\partial\rho\wedge\Psi^F\wedge s\rangle
\end{split}
\end{equation*}
As remarked in the proof of main theorem, $|D_k(\Psi^F\wedge s)|^2=|D_ks|^2+|\bar{\partial}\Psi^F|^2|s|^2$ pointwise, so the only new thing is the computation of boundary terms.

Recall that we can pick $\psi=w_{[I]}$ and then $\Psi^F=|\psi|^{-2}\psi^*\otimes\psi$. Split into two cases.
\begin{enumerate}[(a)]
    \item For $\mathfrak{n}^{1,0}=(\partial\rho)^*\in E_+$, we see $\mathscr{E}_{\partial\rho}\Psi^F=0$ along $\partial M$, and in this case the assumption is $\mathscr{I}_{\partial\rho}s=0$ along $\partial M$. Two boundary terms vanish and the remainging terms are:
    \[\int_{\partial M}\langle\bar{\partial}(\Psi^F\wedge s),\bar{\partial}\rho\wedge \Psi^F\wedge s\rangle+\int_{\partial M}\langle\partial\rho\wedge(\nabla^{1,0*}+\mathcal{T}^*)(\Psi^F\wedge s),\Psi^F\wedge s\rangle\]
    use the formula
    \[\nabla^{1,0*}+\mathcal{T}^*=-\mathscr{I}_{\bar{w}^j}\tilde{\nabla}_{\bar{j}}-T_{\bar{j}\bar{s}r}\mathscr{E}_{\bar{w}^j}\mathscr{I}_{\bar{w}^s}\mathscr{I}_{w^r}\]
    Expand the integrand into convariant differentiations and torsion terms
    \[\langle\bar{\partial}(\Psi^F\wedge s),\bar{\partial}\rho\wedge \Psi^F\wedge s\rangle=\langle\mathscr{E}_{\bar{w}^j}\tilde{\nabla}_{\bar{j}}s,\mathscr{E}_{\bar{\partial}\rho}s\rangle+\frac{1}{2}T_{\bar{j}\bar{s}r}\langle\mathscr{E}_{\bar{w}^j}\mathscr{E}_{\bar{w}^s}\mathscr{I}_{w^r}s,\mathscr{E}_{\bar{\partial}\rho}s\rangle\]
\begin{equation*}
\begin{split}
\langle\partial\rho\wedge(\nabla^{1,0*}+\mathcal{T}^*)(\Psi^F\wedge s),\Psi^F\wedge s\rangle
&=-\langle\mathscr{E}_{\partial\rho}\mathscr{I}_{\bar{w}^j}\tilde{\nabla}_{\bar{j}}\Psi^F,\Psi^F\rangle |s|^2-\langle\mathscr{E}_{\partial\rho}\mathscr{I}_{\bar{w}^j}\Psi^F,\Psi^F\rangle\langle\tilde{\nabla}_{\bar{j}}s,s\rangle\\
&+T_{\bar{j}\bar{s}r}\langle\mathscr{E}_{\partial\rho}\mathscr{I}_{\bar{w}^s}\Psi^F,\Psi^F\rangle\langle\mathscr{E}_{\bar{w}^j}\mathscr{I}_{w^r}s,s\rangle
\end{split}
\end{equation*}
The torsion terms will cancel out:
\[\frac{1}{2}T_{\bar{j}\bar{s}r}\langle\mathscr{E}_{\bar{w}^j}\mathscr{E}_{\bar{w}^s}\mathscr{I}_{w^r}s,\mathscr{E}_{\bar{\partial}\rho}s\rangle=\sum_{i\in[+]}g^{i\bar{j}}\rho_iT_{\bar{j}\bar{s}r}\langle\mathscr{E}_{\bar{w}^s}\mathscr{I}_{w^r}s,s\rangle\]
\[T_{\bar{j}\bar{s}r}\langle\mathscr{E}_{\partial\rho}\mathscr{I}_{\bar{w}^s}\Psi^F,\Psi^F\rangle\langle\mathscr{E}_{\bar{w}^j}\mathscr{I}_{w^r}s,s\rangle=\sum_{i\in[+]}g^{i\bar{s}}\rho_iT_{\bar{j}\bar{s}r}\langle\mathscr{E}_{\bar{w}^j}\mathscr{I}_{w^r}s,s\rangle\]
Next we claim that the derivatives on $s$ will cancel out:
\begin{equation*}
\begin{split}
\langle\mathscr{E}_{\bar{w}^j}\tilde{\nabla}_{\bar{j}}s,\mathscr{E}_{\bar{\partial}\rho}s\rangle-\langle\mathscr{E}_{\partial\rho}\mathscr{I}_{\bar{w}^j}\Psi^F,\Psi^F\rangle\langle\tilde{\nabla}_{\bar{j}}s,s\rangle
&=\sum_{i\in[+]}\rho_i\langle (\mathscr{I}_{w^i}\mathscr{E}_{\bar{w}^j}-g^{i\bar{j}})\tilde{\nabla}_{\bar{j}}s,s\rangle\\
&=-\sum_{i\in[+]}\rho_i\langle\mathscr{E}_{\bar{w}^j}\mathscr{I}_{w^i}\tilde{\nabla}_{\bar{j}}s,s\rangle\\
&=-\langle\mathscr{E}_{\bar{w}^j}\mathscr{I}_{\partial\rho}\tilde{\nabla}_{\bar{j}}s,s\rangle
\end{split}
\end{equation*}
Notice that $\mathscr{I}_{w^j}(s)w_j(\rho)=0$ means that $\mathscr{I}_{w^j}(s)w_j\in T\partial M\otimes\bigwedge\nolimits^{0,q-1}\otimes E\otimes L^k$ along $\partial M$, thus we have
\[0=\langle\tilde{\nabla}_{\bar{w}_j}(\mathscr{I}_{\partial\rho}s),\mathscr{I}_{w^j}s\rangle=\langle\mathscr{I}_{\tilde{\nabla}_{\bar{w}_j}\partial\rho}s,\mathscr{I}_{w^j}s\rangle+\langle\mathscr{I}_{\partial\rho}\tilde{\nabla}_{\bar{w}_j}s,\mathscr{I}_{w^j}s\rangle\]
Finally observe
\[\langle\mathscr{E}_{\bar{w}^j}\mathscr{I}_{\tilde{\nabla}_{\bar{w}_j}\partial\rho}s,s\rangle=\langle(\partial\bar{\partial}\rho)(w_i,\bar{w}_j)\mathscr{E}_{\bar{w}^j}\mathscr{I}_{w^i}s,s\rangle\]
Thus the computation ends up saying that the boundary terms are
\[\int_{\partial M}\langle(\partial\bar{\partial}\rho)(w_i,\bar{w}_j)\mathscr{E}_{\bar{w}^j}\mathscr{I}_{w^i}s,s\rangle-\langle\mathscr{E}_{\partial\rho}\mathscr{I}_{\bar{w}^j}\tilde{\nabla}_{\bar{j}}\Psi^F,\Psi^F\rangle|s|^2\]
\item For $\mathfrak{n}^{1,0}\in E_-$ and $s\in B^{0,q}_-$, we have $\mathscr{I}_{\bar{\partial}\rho}\Psi^F=0$ and $\mathscr{E}_{\bar{\partial}\rho}s=0$ along $\partial M$. This time the boundary integrals are:
\[-\int_{\Gamma_-}\langle\bar{\partial}\rho\wedge(\bar{\partial}^*+\bar{\mathcal{T}}^*)(\Psi^F\wedge s),\Psi^F\wedge s\rangle+\langle\nabla^{1,0}(\Psi^F\wedge s),\partial\rho\wedge\Psi^F\wedge s\rangle\]
Use 
\[\bar{\partial}^{*}+\bar{\mathcal{T}}^*=-\mathscr{I}_{w^i}\tilde{\nabla}_i-T_{ip\bar{q}}\mathscr{E}_{w^i}\mathscr{I}_{w^p}\mathscr{I}_{\bar{w}^q}\]
\[\langle\bar{\partial}\rho\wedge(\bar{\partial}^*+\bar{\mathcal{T}}^*)(\Psi^F\wedge s),\Psi^F\wedge s\rangle=-\langle\mathscr{E}_{\bar{\partial}\rho}\mathscr{I}_{w^i}\tilde{\nabla}_is,s\rangle+T_{ip\bar{q}}\langle\mathscr{E}_{w^i}\mathscr{I}_{\bar{w}^q}\Psi^F,\Psi^F\rangle\langle\mathscr{E}_{\bar{\partial}\rho}\mathscr{I}_{w^p}s,s\rangle\]
\begin{equation*}
\begin{split}
\langle\nabla^{1,0}(\Psi^F\wedge s),\partial\rho\wedge\Psi^F\wedge s\rangle
&=\langle\mathscr{E}_{w^i}\tilde{\nabla}_i\Psi^F,\mathscr{E}_{\partial\rho}\Psi^F\rangle|s|^2+\langle\mathscr{E}_{w^i}\Psi^F,\mathscr{E}_{\partial\rho}\Psi^F\rangle\langle\tilde{\nabla}_is,s\rangle\\
&+\frac{1}{2}T_{ip\bar{q}}\langle\mathscr{E}_{w^i}\mathscr{E}_{w^p}\mathscr{I}_{\bar{w}^q}\Psi^F,\mathscr{E}_{\partial\rho}\Psi^F\rangle|s|^2
\end{split}
\end{equation*}
Torsion terms again cancel out:
\[T_{ip\bar{q}}\langle\mathscr{E}_{w^i}\mathscr{I}_{\bar{w}^q}\Psi^F,\Psi^F\rangle\langle\mathscr{E}_{\bar{\partial}\rho}\mathscr{I}_{w^p}s,s\rangle=\sum_{j\in[-]}g^{p\bar{j}}\rho_{\bar{j}}T_{ip\bar{q}}\langle\mathscr{E}_{w^i}\mathscr{I}_{\bar{w}^q}\Psi^F,\Psi^F\rangle|s|^2\]
\[\frac{1}{2}T_{ip\bar{q}}\langle\mathscr{E}_{w^i}\mathscr{E}_{w^p}\mathscr{I}_{\bar{w}^q}\Psi^F,\mathscr{E}_{\partial\rho}\Psi^F\rangle|s|^2=\sum_{j\in[-]}g^{i\bar{j}}\rho_{\bar{j}}T_{ip\bar{q}}\langle\mathscr{E}_{w^p}\mathscr{I}_{\bar{w}^q}\Psi^F,\Psi^F\rangle|s|^2\]
This time $\mathscr{E}_{\bar{w}^i}(s)w_{\bar{i}}(\rho)=0$ and $\mathscr{E}_{\bar{w}^i}(s)w_{\bar{i}}\in T\partial M\otimes\bigwedge\nolimits^{0,q+1}\otimes E\otimes L^k$ along $\partial M$, so
\[0=\langle\tilde{\nabla}_{w_i}(\mathscr{E}_{\bar{\partial}\rho}s),\mathscr{E}_{\bar{w}^i}s\rangle=\langle\mathscr{E}_{\tilde{\nabla}_{w_i}\bar{\partial}\rho}s,\mathscr{E}_{\bar{w}^i}s\rangle+\langle\mathscr{E}_{\bar{\partial}\rho}\tilde{\nabla}_{w_i}s,\mathscr{E}_{\bar{w}^i}s\rangle\]
\[\langle\mathscr{I}_{w^i}\tilde{\nabla}_{w_i}(\mathscr{E}_{\bar{\partial}\rho}s),s\rangle=\langle(\partial\bar{\partial}\rho)(w_i,\bar{w}_j)\mathscr{I}_{w^i}\mathscr{E}_{\bar{w}^j}s,s\rangle\]
Combining all computations, the final boundary terms are
\[\int_{\partial}\langle(\partial\bar{\partial}\rho)(w_i,\bar{w}_j)\mathscr{I}_{w^i}\mathscr{E}_{\bar{w}^j}s,s\rangle-\langle\mathscr{I}_{\bar{\partial}\rho}\mathscr{E}_{w^i}\tilde{\nabla}_i\Psi^F,\Psi^F\rangle|s|^2\]

\end{enumerate}

\end{proof}

\end{thm}

It turns out that the boundary term has same type of error as in (\ref{thm:Main theorem applied to important bundle}), and the estimate can be obtained if we assume sufficient pseudoconvexity of the boundary.

\begin{example}
Similar to the example (\ref{example:complete product complete}), but we assume that $M_1$ is a weakly pseudoconvex domain with boundary, that is, $\partial\bar{\partial}\rho(v,\bar{v})\ge0$ for $v\in T^{1,0}_x$ and $x\in\partial M$. Since $M_2$ is complete (without boundary), the condition $\mathfrak{n}^{1,0}\in E_+$ along $\partial M$ is satisfied. Under the corresponding domain of $D_k$, we have the boundary term $\ge0$, and again obtain local spectral gap of growth order $k$ (\ref{local spectral gap assumption}) for $U\subset\subset M(q)$ if for $k>>0$,
\[\prescript{}{k}{f}_{d-q,q+1}=\frac{2}{3}k\lambda_{d-q}-\frac{2}{3}\lambda^E_{d-q}-C_d|\partial\Theta|^2\ge0\qquad\prescript{}{k}{f}_{d-q,q-1}=-\frac{2}{3}k\lambda_{d-q+1}-\frac{2}{3}\lambda^E_{d-q}-C_d|\partial\Theta|^2\ge0\]
In particular if $M_1$ is a weakly pseudoconvex domain in $\C^{d-q}$ and $M_2=\C^q$ and $L^k$-valued forms, always
\[\prescript{}{k}{f}_{d-q,q+1}=k\lambda_{d-q}\ge0\qquad\prescript{}{k}{f}_{d-q,q-1}=-k\lambda_{d-q+1}\ge0\]
Similar conclusion holds for $M_1$ complete and $M_2$ weakly pseudoconvex domain.
\end{example}

\subsection{Full asymptotic expansion}

Notice that all cases we discuss before satisfies $\ker D_k=\ker\Box_k$, and if the geometric conditions are satisfied then we have local spectral gap of growth order $k$ (\ref{local spectral gap assumption}) on any $U\subset\subset M(q)$ since $\prescript{}{k}{f}_{p,q}\ge C_Uk$ on $U$ for some $C_U>0$.

For any $z_0\in M$, we will always take a local coordinate as follows:
\begin{itemize}
    \item $z_0$ is identified with $0\in\C^d$.
    \item $g_{i\bar{j}}(z)=g_0+O(|z|)$, $(g_0)_{i\bar{j}}=\delta_{ij}$.
    \item $L$ is trivialized by section $s$ with metric represented by the positive function $e^{-\phi}=|s|^2$, and $\phi(z)=\phi_0(z)+R_3(z)$, $\phi_0=\sum_{i=1}^d\lambda_i(z_0)|z^i|^2$ and $R_3(z)=O(|z|^3)$.
    \item $E$ is trivialized by local holomorphic frame $\{e_a\}_{a=1}^r$ with metric $h(z)=h_0+O(|z|)$, $(h_0)_{a\bar{b}}=\delta_{ab}$.
    
\end{itemize}

Let us begin with an easy lemma.
\begin{lemma}\label{skullduggery}
Let $\tilde{L},\tilde{E}$ be the trivial bundles of rank $1$ and $r$ over $\C^d$.

For any $z_0\in M(q)$, there exists an $\epsilon>0$ and metrics $\tilde{g},\tilde{\phi},\tilde{h}$ on $\C^d,\tilde{L},\tilde{E}$ so that
\[(\tilde{g}|_{B_{3\epsilon}(0)},\tilde{\phi}|_{B_{3\epsilon}(0)},\tilde{h}|_{B_{3\epsilon}(0)})=(g^M|_{B_{3\epsilon}(0)},\phi^L|_{B_{3\epsilon}(0)},h^E|_{B_{3\epsilon}(0)})\]
under the identification with $B_{3\epsilon}(0)\subset\C^d$, and that outside $B_{4\epsilon}(0)$ we have
\[g_{i\bar{j}}(z)=g_0,\ \tilde{\phi}(z)=\sum_{i=1}^d\lambda_i(z_0)|z^i|^2,\ \tilde{h}(z)=h_0\]
Moreover, the induced curvature $R^{\tilde{L}}(z)$ is nondegenerate with signature $(n-q,q)$ for every $z\in\C^d$.
\begin{proof}
The only problem is to make $R^{\tilde{L}}$ nondegenerate with fixed signature.

There exists constant $C_1,C_2>0$ such that for any $\epsilon>0$ small we can take a smooth cut off $\zeta:\C^d\to[0,1]$ with $\zeta=1$ on $B_{3\epsilon}(0)$, $\zeta=0$ outside $B_{4\epsilon}(0)$ and that
\[|d\zeta|\le\frac{C_1}{\epsilon}\qquad|\partial\bar{\partial}\zeta|\le\frac{C_2}{\epsilon^2}\]
Define $\tilde{g}=\zeta g+(1-\zeta)g_0$ and $\tilde{h}=\zeta h+(1-\zeta)h_0$ works.

Consider the function $\tilde{\phi}=\zeta\phi+(1-\zeta)\phi_0=\phi_0+\zeta\cdot R_3$ defined on $\C^d$, we want to take $\epsilon$ small enough so that $\partial\bar{\partial}\tilde{\phi}$ also has same signature with $\partial\bar{\partial}\phi_0$. But notice that
\[\partial\bar{\partial}(\zeta\cdot R_3)=R_3\partial\bar{\partial}\zeta+\zeta\partial\bar{\partial}R_3+\partial\zeta\wedge\bar{\partial}R_3+\partial R_3\wedge\bar{\partial}\zeta\]
and the norm is bounded by $C\cdot\epsilon$ for $|z|\le 4\epsilon$. Thus such $\epsilon$ exists.
    
\end{proof}

\end{lemma}

\begin{thm}[Localization]\label{thm:localization}
Assume (\ref{local spectral gap assumption}). For any $z_0\in U$, there exists an $\epsilon>0$ such that:
Given any $l,N\in\N$, there exists a $C_{l,N}>0$ such that
\[|F_{\epsilon}(D_k)(z,z^\prime)-B^{(q)}_k(z,z^\prime)|_{C^l(B_{\epsilon}(z_0)\times B_{\epsilon}(z_0))}\le C_{l,N}k^{-N}\]
where
\[F_{\epsilon}(x)=\int_{-\infty}^{\infty}\chi_{\epsilon}(v)e^{\sqrt{-1}vx/\sqrt{2}}\,dv\]
and $\chi_{\epsilon}\in C^{\infty}_c(\R)$ is an even function supported in $[-\epsilon,\epsilon]$ and $\int_{-\infty}^\infty\chi_{\epsilon}(v)\,dv=1$.

\begin{proof}
Take the $\epsilon>0$ and $(\tilde{g},\tilde{\phi},\tilde{h})$ as in (\ref{skullduggery}).

For any small $a>0$ and $\eta(0,x)\in\Omega^{0,\bullet}_c(B_{a}(0),E\otimes L^k)$, notice that $\eta(t,x)=e^{\sqrt{-1}tD_k/\sqrt{2}}\eta(0,x)$ satisfies
\[(\frac{\partial^2}{\partial t^2}+\Box_k)\eta(t,x)=0\]
with $(\frac{\partial}{\partial t}\omega)(0,x)=\frac{\sqrt{-1}}{\sqrt{2}}D_k\eta(0,x)$ supported in $B_{a}(0)$. By finite propagation speed for wave equation \hyperlink{ref:MM1}{[MM1, Thm.D.2.1]}, we know that $\eta(t,x)$ is supported in $B_{a+t}(0)$, thus
\[F_{\epsilon}(D_k):\Omega^{0,\bullet}_c(B_{a}(0),E\otimes L^k)\to\Omega^{0,\bullet}_c(B_{a+\epsilon}(0),E\otimes L^k)\]
Notice that this operator only depends on $D_k|_{B_{a+\epsilon}(0)}$. Now since $F_{\epsilon}(0)=1$, apply local spectral gap, for $s\in\Omega_c^{0,q}(B_{\epsilon}(0),E\otimes L^k)$,
\[||F_{\epsilon}(D_k)s-B_ks||=||F_{\epsilon}(D_k)s-B_kF_{\epsilon}(D_k)s||\le C_Uk^\alpha||\Box_k F_{\epsilon}(D_k)s||\]
Notice that \hyperlink{ref:MM1}{[MM1, Lem.1.6.2]} we have for $s\in H^{2l}(B_{\epsilon}(0),E\otimes\bigwedge\nolimits^{0,q}\otimes L^k)$,
\[||s||_{H^{2l}(B_{\epsilon}(0))}\le C\sum_{j=0}^lk^{4(l-j)}||\Box_k^js||_{L^2(B_{2\epsilon}(0))}\]
Thus for differential operators, $P$ with degree smaller that $2l$, $Q$ with scalar principal symbol, and both compactly supported in $B_{\epsilon}(0)$,
\begin{equation*}
\begin{split}
||P(F_{\epsilon}(D_k)-B_k)Qs||_{L^2(B_{\epsilon}(0))}
&\le C(k^{4l}||(F_{\epsilon}(D_k)-B_k)Qs||_{L^2(B_{2\epsilon}(0))}+\sum_{j=1}^{l}k^{4(l-j)}||\Box_{k}^j(F_{\epsilon}(D_k)-B_k)Qs||_{L^2(B_{2\epsilon}(0))})\\
&\le C(C_Uk^{4l+\alpha}||\Box_k F_{\epsilon}(D_k)Qs||_{L^2(B_{2\epsilon}(0))}+\sum_{j=1}^lk^{4(l-j)}||\Box^{j}_kF_{\epsilon}(D_{k})Qs||_{L^2(B_{2\epsilon}(0))})
\end{split}
\end{equation*}
Now $\Box_k^jF_{\epsilon}(D_k)Qs$ is supported in $B_{2\epsilon}(0)$ and for $s^\prime\in\Omega^{0,\bullet}_c(B_{2\epsilon}(0),E\otimes L^k)$,
\[\langle\Box^j_k F_{\epsilon}(D_k)Qs,s^\prime\rangle=\langle s, Q^* F_{\epsilon}(D_k)\Box^j_k s^\prime\rangle\]
Using again \hyperlink{ref:MM1}{[MM1, Lem.1.6.2]} and local spectral gap condition, we only need to estimate $||\Box_k^mF_{\epsilon}(D_{k})\Box_k^js||$ for $s\in\Omega^{0,q}_c(B_{2\epsilon}(0),E\otimes L^k)$ and $m>0$. Notice that the operator $\Box_k^mF_{\epsilon}(D_{k})\Box_k^j$ acting on $\Omega^{0,\bullet}_c(B_{2\epsilon}(0),E\otimes L^k)$ only depends on the restriction of $D_k$ to $B_{3\epsilon}(0)$, we know that is is equal to $\Box_{k\tilde{\phi}}^mF_{\epsilon}(D_{k\tilde{\phi}})\Box_{k\tilde{\phi}}^j$. But by the estimate for complete case (\ref{model case}), $\Spec D_{k\tilde{\phi}}\subset\{0\}\cup[C\sqrt{k},\infty)$ for some $C>0$, and $F_{\epsilon}$ is a Schwartz function, so for any $N\in\N$, there exists $C_N>0$ such that
\[||\Box_{k\tilde{\phi}}^mF_{\epsilon}(D_{k\tilde{\phi}})\Box_{k\tilde{\phi}}^js||=||\Box_{k\tilde{\phi}}^m(F_{\epsilon}(D_{k\tilde{\phi}})-B_{k\tilde{\phi}})\Box_{k\tilde{\phi}}^js||\le C_Nk^{-N}||s||\]
Combining all and we get the result.

\end{proof}

\end{thm}

Notice that apply the same localization to the metrics $(\tilde{g},\tilde{\phi},\tilde{h})$ on $(\C^d,\tilde{L},\tilde{E})$, we know that the asymptotic behavior of $B_k^{(q)}(z,z^\prime)$ near $(z_0,z_0)$ is equivalent to the Bergman kernel $B_{k\tilde{\phi}}^{(q)}(z,z^\prime)$ of the model case (\ref{model case}). In this case we have a complete manifold with $\Spec D_k\subset\{0\}\cup[C\sqrt{k},\infty)$. $D_k$ and $\Box_k$ is defined in Gaffney sense, and we can apply the result in \hyperlink{ref:HM}{[HM]} to obtain a full diagonal expansion.

Since we have global spectral gap, we also can use the method in \hyperlink{ref:MM1}{[MM1]} to obtain the expansion. We briefly explain how the method applies. They consider the modified Dirac operator $D^{c,A_0}_k$ which coincides with $D_k$ near $0$, and outside a compact set we have $(\tilde{g},\tilde{\phi},\tilde{h})=(g_0,\phi_0,h_0)$, thus we also have a spectral gap of growth order $k$ for $(D^{c,A_0}_k)^2$. Same localization works and we only need to investigate the kernel of orthogonal projection onto $\ker(D^{c,A_0}_k)^2$. Put $t=\frac{1}{\sqrt{k}}$ and define the appropriate rescaled operator $\mathscr{L}^t_2$ for $(D^{c,A_0}_k)^2$. Now $\Spec\mathscr{L}^t_2\subset\{0\}\cup[\mu_0,\infty)$ and we have for any $l\ge1$, the 
orthogonal projection $P_{0,t}$ to $\ker\mathscr{L}^t_2$ can be written as
\[P_{0,t}=\frac{1}{2\pi\sqrt{-1}}\int_{\delta}\lambda^{l-1}(\lambda-\mathscr{L}^t_2)^{-l}\,d\lambda\]
where $\delta$ is the circle of radius $\mu_0/4$ oriented counterclockwise. Use Taylor expansion to the coefficient of $\mathscr{L}^t_2$ and write
\[\mathscr{L}^t_2=\mathscr{L}^0_2+\sum_{r=1}^mt^r\mathcal{O}_r+O(t^{m+1})\]
Do a sequence of Sobolev estimate on the resolvent $(\lambda-\mathscr{L}^t_2)^{-1}$ so that we can insert the expansion uniformly into the residue formula, we can get:

\begin{thm}[\hyperlink{ref:MM1}{[MM1, Thm.8.2.4]}]
Assume (\ref{local spectral gap assumption}), there exists smooth coefficients $b_r(x)\in\operatorname{End}(\bigwedge\nolimits^{0,\bullet}\otimes E)_x$ 
with the following property: for any $m,n\in\N$ and $U^{\prime}\subset\subset U$, there exists $C_{m,n,U^\prime}>0$ such that for $k>>0$,
\[|B_k(x,x)-\sum_{r=0}^mb_r(x)k^{d-r}|_{C^n(U^\prime)}\le C_{m,n,U^\prime}k^{d-m-1}\]
In particular,
\[b_0=|\det(R^L/(2\pi))|I_{\det(\overline{W}^*)\otimes E}\]
where $W$ is the subbundle of $T^{1,0}$ spanned by the eigenvectors of the negative eigenvalues of $R^L\in\operatorname{End}(T^{1,0})$.
\end{thm}
We also have off-diagonal expansion. For more information about the coefficients and details, see \hyperlink{ref:MM1}{[MM1, Thm.4.1.7-21]}.

\section*{References}

\hypertarget{ref:Aron}{[Aron]} N. Aronszajn, \textit{A unique continuation theorem for solutions of elliptic partial differential equations or inequalities of second order}, J. Math. Pure. Appl., 36 (1957), 235-249.

\hypertarget{ref:AMM}{[AMM]} H. Auvray, X. Ma and G. Marinescu, \textit{Bergman kernels on punctured Riemann surfaces}, Math. Ann.
379 (2021), 951–1002.

\hypertarget{ref:Berm}{[Berm]} R. Berman, \textit{Bergman kernels and local holomorphic Morse inequalities}, Math. Z. 248 (2004), no. 2,
325–344.

\hypertarget{ref:Bern}{[Bern]} B. Berndtsson, \textit{$L^2$-methods for the $\bar{\partial}$-equation}, KASS University Press-Johanneberg-Masthugget-Sisjön, 1995.

\hypertarget{ref:Dav}{[Dav]} E.-B. Davies, \textit{Spectral Theory and Differential Operators}, Cambridge Stud. Adv. Math., vol. 42, (1995).

\hypertarget{ref:DLM}{[DLM]} X. Dai, K. Liu and X. Ma, \textit{On the asymptotic expansion of Bergman kernel}, C. R. Math. Acad. Sci. Paris
339 (2004), no. 3, 193–198.

\hypertarget{ref:GH}{[GH]} P. Griffiths and J. Harris, \textit{Principles of Algebraic Geometry}, John Wiley and
Sons, New York, 1978.

\hypertarget{ref:HH}{[HH]} C.-Y. Hsiao and R.-T. Huang, \textit{G-invariant Szegő kernel asymptotics and CR reduction}, arXiv: 1702.05012.

\hypertarget{ref:HM}{[HM]} C-Y. Hsiao and G. Marinescu, \textit{Asymptotics of spectral function of lower energy forms and
Bergman kernel of semi-positive and big line bundles}, Comm. Anal. Geom. 22(1) (2014), 1-108.

\hypertarget{ref:HMM}{[HMM]} C.-Y Hsiao, X. Ma, and G. Marinescu, \textit{Geometric quantization on CR manifolds}, Commun. Contemp. Math. (2023).

\hypertarget{ref:Hör}{[Hör]} L. Hörmander, \textit{$L^2$ estimates and existence theorems for the $\bar{\partial}$ operator}. Acta Math., 113
(1965), 89–152.

\hypertarget{ref:MM1}{[MM1]} X. Ma and G. Marinescu, \textit{Holomorphic Morse inequalities and Bergman kernels}, Progress in Math., vol.
254, Birkh¨auser, Basel, 2007, 422 pp.

\hypertarget{ref:MM2}{[MM2]} X. Ma and G. Marinescu, \textit{The first coefficients of the asymptotic expansion of the Bergman kernel of the
spin$^c$ Dirac operator}, Internat. J. Math., 17 (2006), no. 6, 737–759.

\hypertarget{ref:MZ}{[MZ]} X. Ma and W. Zhang, \textit{Bergman kernels and symplectic reduction}, Ast´erisque 318 (2008), 154 pp .

\hypertarget{ref:Tian}{[Tian]} G. Tian, \textit{On a set of polarized Kähler metrics on algebraic manifolds}, J. Differential Geom. 32 (1990), 99–130.

\end{document}